\documentclass[11pt]{article}
\usepackage{amsmath}
\usepackage{subfigure}
\usepackage{amsthm}
\usepackage{amstext}
\usepackage{amsopn}
\usepackage{texdraw}
\usepackage{pstricks,pst-plot}
\usepackage{mathrsfs}
\usepackage[numbers,sort&compress]{natbib}
\usepackage{multicol,graphicx,color}
\usepackage{pslatex}
\usepackage{amssymb}
\usepackage{latexsym}
\usepackage{lscape}
\usepackage{epsfig}
\usepackage{amsfonts}
\usepackage{enumerate}
\usepackage{exscale}
\usepackage{relsize}

\numberwithin{equation}{section}
\theoremstyle{plain}
\newtheorem{theorem}{Theorem}[section]
\newtheorem{lemma}[theorem]{Lemma}
\newtheorem{proposition}[theorem]{Proposition}
\newtheorem{corollary}[theorem]{Corollary}
\theoremstyle{definition}

\theoremstyle{remark}
\newtheorem{remark}{Remark}

\begin{document}

\title{Geometric properties of minimizers in the planar three-body problem with two equal masses}

\author{
Wentian Kuang\\
Department of Mathematics, Southern University of Science and Technology\\
Shenzhen 518055, P.R. China\\
Email: kuangwt@sustc.edu.cn
 \and
Duokui Yan \\
School of Mathematical Sciences, Beihang University\\
Beijing 100191, P.R. China \\
Email: duokuiyan@buaa.edu.cn
}
\date{}


\maketitle

\begin{abstract}
It it shown that each lobe of the figure-eight orbit is star-shaped, which implies the polar angle is monotone in each lobe. In general, it is not clear when a minimizer is star-shaped. In this paper, we study minimizers connecting two fixed-ends (i.e. the Bolza problem) in the planar three-body problem with two equal masses. We show that if the Jacobi coordinates of the two fixed-ends are in adjacent closed quadrants, then the corresponding minimizer must stay in two adjacent closed quadrants. If we further assume the two Jacobi coordinates are orthogonal on one of the fixed-ends, then the polar angles of the Jacobi coordinates in the minimizer have at most one critical point. If the two Jacobi coordinates are orthogonal on both ends, then the two polar angles must be monotone. These geometric properties can be applied to show the existence of two sets of periodic orbits.

\end{abstract}


{\bf AMS classification number:} 37N05, 70F10, 70F15, 37N30, 70H05, 70F17\\


\section{Introduction}
Variational method is an important tool in studying periodic orbits in the $N$-body problem. After the pioneering work of the figure-eight orbit \cite{CM} in the equal mass three-body problem, many new periodic orbits have been shown to exist by variational methods. Novel ideas and new theories have been innovated in this direction during the last few decades.

In \cite{CH, CH1}, Chen studied the existence of retrograde and prograde orbits of three different masses in the planar three-body problem. Intuitively, `` a retrograde orbit of the planar three-body problem is a relative periodic (periodic in a rotating coordinate system) solution with two adjacent masses revolving around each other in one direction while their mass center revolves around the third mass in the other direction. An orbit is said to be prograde if both revolutions follow the same direction." \cite{CH1} Rigorous definitions are given by the homology class of the projected curve on the shape sphere in \cite{CH} and by braids in \cite{CH1}. 

It is shown in \cite{CM} that the figure-eight orbit can be characterized as a minimizer connecting an Euler configuration and an isosceles configuration. By introducing the shape sphere coordinates \cite{M1, M2,M4, M3, MOE, MOV}, they prove that this minimizer contains no collinear or isosceles configuration except for the boundaries. The polar angle of each body in this minimizer is monotone.

Meanwhile, from the numerical pictures of the retrograde and prograde orbits in \cite{CH1}, it seems that one of the three bodies has a star-shaped trajectory, which implies that its polar angle could be monotone. Furthermore, if one considers the retrograde orbits with two equal masses, numerically it can be characterized as a minimizer connecting a collinear configuration and an isosceles configuration. It seems that such a minimizer contains no collinear or isosceles configuration except for the boundaries. Motivated by these numerical observations, we study properties of minimizers connecting two fixed-ends (i.e. the Bolza problem) in the planar three-body problem with two equal masses. In fact, we can show that if both fixed ends are isosceles, the polar angles of the Jacobi coordinates in the corresponding minimizer must be monotone.


Let $M=(m_1, \, m_2,\, m_3)$ be the mass vector. We set $q_i\in\mathbb{C}$ to be the position of mass $m_i \in \mathbb{R}^+\, (i=1,2,3)$ and $q=(q_1, \, q_2, \, q_3)$. Without loss of generality, we assume the center of mass to be at the origin. That is, $q \in \Sigma$:
\begin{equation}\label{configuration}
\Sigma=\left\{ q \in \mathbb{C}^{3} \, \bigg{|} \, \sum_{i=1}^3 m_iq_i =0 \right \}.
\end{equation}

The Jacobi coordinates $Z=(Z_1, \, Z_2)\in\mathbb{C}^2$ (as in Fig. \ref{Jacobifigure}) for the planar three-body problem are given by
\begin{equation}\label{jacob01}
Z_1=q_3-q_2,  \qquad \quad Z_2=q_1-\frac{m_2q_2+m_3q_3}{m_2+m_3}=(1+\frac{m_1}{m_2+m_3})q_1.
\end{equation}

\begin{figure}[ht]
	\begin{center}
		\psset{xunit=2.4cm,yunit=2.4cm}
		\begin{pspicture}(-1,-1)(1,1)
		\psline[linewidth=1pt, linecolor=blue](-0.3,0.3)(-0.07, 0.07)
		\psline[linewidth=1pt, linecolor=blue](-0.8,0.15)(0.2,0.45)
		\psdots[dotsize=4pt](0.4,-0.4)(-0.8,0.15)(0.2,0.45)(-0.3,0.3)(-0.07, 0.07)
		\psaxes[linewidth=1.5pt]{->}(-0.07, 0.07)(-0.95, -0.85)(0.95, 0.85)
		\psline[linewidth=1pt, linecolor=blue, linestyle=dashed]{->}(-0.07, 0.07)(0.93, 0.37)
		\psline[linewidth=1pt, linecolor=blue, linestyle=dashed]{->}(-0.07, 0.07)(0.63, -0.63)
		\rput(1, -0.1){$x$}
		\rput(0.1, 0.85){$y$}
		\rput(0.98, 0.47){\textcolor{blue}{$Z_1$}}
		\rput(0.5, -0.73){\textcolor{blue}{$Z_2$}}
		\rput(-0.17, -0.1){0}
		\rput(-0.95, 0.25){$q_2$}
		\rput(0.38, 0.55){$q_3$}
		\rput(0.55, -0.3){$q_1$}
		\end{pspicture}
	\end{center}
	\caption{ \label{Jacobifigure}The Jacobi coordinates, where $Z_1$, $Z_2$ are vectors in dashed lines. }
\end{figure}


Let $Q_0$ and $Q_1$ be two given configurations in $\Sigma$. By \cite{CA, Mar}, there exists a minimizer $\mathcal{P}_{Q_0 Q_1} \in H^1([0,\,1], \, \Sigma)$, such that
\begin{equation}\label{Bolza}
\mathcal{A}(\mathcal{P}_{Q_0 Q_1} )= \inf_{q \in  P(Q_0,\, Q_1)} \mathcal{A}(q) = \inf_{q \in  P(Q_0, \, Q_1)} \int_0^1 (K+U) \, dt,
\end{equation}
where $\displaystyle K=\frac{1}{2}\sum_{i=1}^{3}m_i |\dot{q}_i|^2$ is the kinetic energy, $\displaystyle U={\sum_{\substack{1\le i< j\le 3}}\frac{m_im_j}{|q_i-q_j|}}$ is the potential function, and
\[ \mathrm P(Q_0,\, Q_1)=\left\{q\in H^1([0,\, 1],\, \Sigma) \, \big{|}  \, q(0)=Q_0, \,\, q(1) =Q_1 \right\}.\]
It is known that the minimizer $\mathcal{P}_{Q_0 Q_1}$ is a weak solution of the Newtonian equations:
\begin{equation}\label{newton}
m_i\ddot{q}_i=\frac{\partial U}{\partial q_i},\quad \quad \ (i=1, 2, 3),
\end{equation}
where $\ddot{q}_i=\frac{d^2 q_i}{dt^2}$ and $U$ is the potential function as in \eqref{Bolza}.

In what follows, we assume the masses are $M=(1, \, m, \, m)$. The Jacobi coordinates \eqref{jacob01} becomes
\[ Z_1=q_3-q_2,  \qquad \quad Z_2=q_1-\frac{q_2+q_3}{2}=(1+\frac{1}{2m})q_1. \]
Note that the map from $q= (q_1, \, q_2,\, q_3)\in \Sigma$ to $Z=(Z_1,Z_2)\in \mathbb{C}^2$ is linear and it is bijective. Therefore, we can analyze the action functional \eqref{Bolza} and the Newtonian equations \eqref{newton} under the Jacobi coordinates.

In fact, the kinetic energy and the potential function can be rewritten as
\begin{equation}\label{Kformulajacob}
K=\frac{1}{2} \sum_{i=1}^3 m_i |\dot{q}_i|^2= \frac{m}{4} |\dot{Z}_1|^2+ \frac{m}{2m+1} |\dot{Z}_2|^2,
\end{equation}
\begin{equation}\label{Uformulajacob}
U= \sum_{1 \le i<j\le 3} \frac{m_i m_j}{|q_i-q_j|}=\frac{m^2}{|Z_1|}+ \frac{m}{|\frac{1}{2}Z_1+ Z_2|}+ \frac{m}{|\frac{1}{2}Z_1-Z_2|}.
\end{equation}
And the Newtonian equations \eqref{newton} under the Jacobi coordinates are
\begin{equation}\label{Z1Z2equation0}
\begin{split}
&\ddot{Z}_1=\frac{Z_2-Z_1/2}{|Z_2-Z_1/2|^3} - \frac{Z_2+Z_1/2}{|Z_2+Z_1/2|^3} - \frac{2mZ_1}{|Z_1|^3},\\
&\ddot{Z}_2=-\frac{1+2m}{2}\left[\frac{Z_2-Z_1/2}{|Z_2-Z_1/2|^3}+\frac{Z_2+Z_1/2}{|Z_2+Z_1/2|^3}\right].
\end{split}
\end{equation}

We then define some notations and introduce the main results.

\textbf{Notations:}
\begin{itemize}
	\item $\mathsf{Q}_i$ is the $i$-th $(i=1, 2, 3, 4)$ quadrant in the Cartesian $xy$ coordinate system and $\overline{\mathsf{Q}_i}$ is its closure. For example, $\mathsf{Q}_1=\{ (x, \,  y) \, |\, x>0,\,  y>0 \}$ and $\overline{\mathsf{Q}_1}=\{(x, \, y) \, | \, x\geq 0, \, y \geq 0\}$.
	 \item $(r_j, \, \theta_j)$ is the polar coordinates of $Z_j$: \, $Z_j=r_j e^{ i\, \theta_j}$ with $r_j \ge 0  \, (j=1,2)$.
     \item  $Z_1(t_0) \perp Z_2(t_0)$ if the real part of $Z_1(t_0)  \overline{Z_2(t_0) }$ is $0$, i.e. $Re \left(Z_1(t_0)  \overline{Z_2(t_0) } \right)=0$.
     \item A minimizer $Z=(Z_1, \, Z_2) \, (t \in [0,\, 1])$ is called \textbf{an isosceles triangular motion} if $Z_1(t) \perp Z_2(t)$ and $\dot{\theta}_1(t)=\dot{\theta}_2(t)\equiv 0$ for all $t \in (0,1)$. It is said to be \textbf{a collinear motion} if $Z_1(t) \parallel Z_2(t)$ and $Z_2(t)  \not\equiv 0$ for all $t \in (0,\,1)$.
	\item For $a\le b$,
	\begin{equation}\label{cone}
	\begin{split}
	C_{[a,\, b]}&=\left\{Z=re^{i\theta}\in \mathbb{C}\, \, \Big| \, \,  r \ge 0, \, \, \theta\in [a, \, b] \right\};\\
	C_{[a,\, b]}^+&=C_{[a,\, b]}\cdot e^{\frac{i\pi}{2}}=\left\{Z=re^{i\theta}\in \mathbb{C}\, \, \Big| \, \,r \ge 0,  \, \, \theta\in \left[a+\frac{\pi}{2}, \, b+\frac{\pi}{2} \right]\right\};\\
	C_{[a,\, b]}^-&=C_{[a,\, b]}\cdot e^{-\frac{i\pi}{2}}=\left\{Z=re^{i\theta}\in \mathbb{C}\, \, \Big| \, \, r \ge 0, \, \, \theta\in \left[a-\frac{\pi}{2}, \, b-\frac{\pi}{2} \right] \right\}.
	\end{split}
	\end{equation}
\end{itemize}

The main results are as follows.

\begin{theorem}[Theorem \ref{thmquadrants} in Section 2]\label{maingeometricresult1}
Assume $Z_1(0), \, Z_1(1)\in \overline{\mathsf{Q}_i} \, $ and $ \, Z_2(0), \, Z_2(1)\in \overline{\mathsf{Q}_j}$, where $\overline{\mathsf{Q}_i}$ and $\overline{\mathsf{Q}_j}$ are two adjacent closed quadrants. Let $Z=(Z_1, \,Z_2) \in H^1([0, \,1], \mathbb{C}^2)$ be a minimizer connecting the two fixed-ends: $Z(0)=(Z_1(0), \, Z_2(0))$ and $Z(1)=(Z_1(1), \, Z_2(1))$. That is,
\begin{equation}\label{Bolza2}
 \mathcal{A}(Z)=\inf_{\widehat{Z} \in P(Z(0), \, Z(1))}\mathcal{A}(\widehat{Z})= \inf_{\widehat{Z} \in P(Z(0), \, Z(1))} \int_0^1 (K+U) \, dt,
 \end{equation}
where $K$ (in \eqref{Kformulajacob}) and $U$ (in \eqref{Uformulajacob}) are the kinetic energy and potential function respectively, and
\[P(Z(0), \, Z(1))= \left\{\, \widehat{Z}=(\widehat{Z}_1, \, \widehat{Z}_2) \in H^1([0,1], \mathbb{C}^2) \, \big| \, \widehat{Z}(0)=Z(0), \, \,  \widehat{Z}(1)=Z(1) \, \right\}. \]
Then $Z_1(t)$ and $Z_2(t)$ must be in two adjacent closed quadrants for all $t \in [0,\, 1]$ and they must satisfy one of the following three cases:\\
\hspace*{0.2in} (a) both $Z_1(t)$ and $Z_2(t)$ are away from the coordinate axes for all $t \in (0, \, 1)$; \\
\hspace*{0.2in}  (b) both $Z_1(t)$ and $Z_2(t)$ are on the coordinate axes for all $t\in [0, \,1]$ and they never touch the origin in $(0,\, 1)$;  \\
\hspace*{0.2in}  (c) $Z=(Z_1, \, Z_2)$ is part of an Euler orbit with $Z_2(t) \equiv 0$ for all $t \in [0,\,  1]$.
\end{theorem}

\begin{figure}[ht]
\begin{center}
\psset{xunit=2.4cm,yunit=2.4cm}
\begin{pspicture}(-2.4,-0.3)(2.4,1)
\psaxes[linewidth=1.5pt]{->}(-1.4, 0.07)(-2.2, -0.25)(-0.5, 0.9)
\psaxes[linewidth=1.5pt]{->}(1.4, 0.07)(0.6, -0.25)(2.3, 0.9)
\psline[linecolor=cyan, doubleline=true]{->}(-0.15, 0.07)(0.15, 0.07)
\psdots[dotsize=4pt](-0.8, 0.3)(-1.3, 0.7)(-1.7, 0.2)(-1.9, 0.8)
\rput(-0.6, 0.4){\textcolor{blue}{\small{$Z_1(0)$}}}
\rput(-1.1, 0.8){\textcolor{blue}{\small{$Z_1(1)$}}}
\rput(-1.95, 0.25){\textcolor{red}{\small{$Z_2(0)$}}}
\rput(-2.1, 0.93){\textcolor{red}{\small{$Z_2(1)$}}}
\psdots[dotsize=4pt](2, 0.3)(1.5, 0.7)(1.1, 0.2)(0.9, 0.8)
\pscurve[linewidth=2pt, linecolor=red, linestyle=dashed]{->}(1.1, 0.2)(1, 0.25)(0.7, 0.3)(0.9, 0.8)
\pscurve[linewidth=2pt, linecolor=blue, linestyle=dotted]{->}(2, 0.3)(1.9, 0.8)(1.5, 0.7)
\rput(0.5, 0.6){\textcolor{red}{\small{$Z_2(t)$}}}
\rput(2.25, 0.6){\textcolor{blue}{\small{$Z_1(t)$}}}
\rput(-1.3, -0.05){0}
\rput(1.5, -0.05){0}
 \end{pspicture}
   \end{center}
 \caption{ \label{Jacobifigure2}An illustration of case (a) in Theorem  \ref{maingeometricresult1}. On the left graph, we assume $Z_1(0), Z_1(1) \in \mathsf{\overline{Q}}_1$ and  $Z_2(0), Z_2(1) \in \mathsf{\overline{Q}}_2$. Theorem  \ref{maingeometricresult1} implies that $Z_1(t)\in \mathsf{\overline{Q}}_1$ and $Z_2(t)\in \mathsf{\overline{Q}}_2$ for all $t \in [0, \,1]$. }
\end{figure}
Fig. \ref{Jacobifigure2} illustrates case (a) in the above theorem. In fact, part of Theorem \ref{maingeometricresult1} has been shown in \cite{Yan3}. Here we give a complete proof by considering equations \eqref{Z1Z2equation0} under the polar coordinates.

The next corollary can be seen as an extension of Theorem \ref{maingeometricresult1}.
\begin{corollary}[Corollary \ref{thmcone} in Section 2]\label{corcone}
		Assume there exist $a, b \in \mathbb{R}$ satisfying $0\le b-a<\frac{\pi}{2}$, such that $Z_1(0), \,Z_1(1)\in C_{[a,\,b]}$, $Z_2(0),\, Z_2(1)\in C_{[a,\,b]}^+$. We further assume that $|Z_i(0)|+|Z_i(1)| \ne 0, \, (i=1,\,2)$. Let $Z=(Z_1, \,Z_2) \in H^1([0, \,1], \mathbb{C}^2)$ be a minimizer connecting the two fixed-ends: $Z(0)$ and $Z(1)$.

Then $Z_1(t)\in C_{[a,\,b]}$ and $Z_2(t)\in C_{[a,\,b]}^+$ for all $t\in [0,\,1]$, and they must satisfy one of the following two cases:\\
\hspace*{0.4in}	    (a) both $Z_1(t)$ and $Z_2(t)$ are away from the boundaries of the corresponding closed cones for all $ t \in (0, \,1)$;\\
\hspace*{0.4in} 	(b) both $Z_1(t)$ and $Z_2(t)$ are on the boundaries of the corresponding closed cones for all $t\in [0, \,1]$, and they can not reach the origin when $ t \in (0, \,1)$.

Similar conclusion holds if $Z_2(0), \, Z_2(1)\in C_{[a,\,b]}^-$.
\end{corollary}
\begin{remark}
In Corollary \ref{corcone}, we impose an assumption that $|Z_i(0)|+|Z_i(1)| \ne 0, \, (i=1,\,2)$. Consequently, the closed cones of $Z_1$ and $Z_2$ can be determined by the boundary conditions. It is clear that case $(c)$ in Theorem \ref{maingeometricresult1} is eliminated by the assumption. Indeed, if $Z_i(0)=Z_i(1)=0$ for $i=1$ or $2$, we can still prove that there exist two orthogonal closed cones so that $Z_1(t)$ belongs to one cone and $Z_2(t)$ belongs to the other for all $t\in [0,1]$, but the two cones may not be $C_{[a,\,b]}$ and $C_{[a,\,b]}^+$.
\end{remark}

Under the assumptions in Theorem \ref{maingeometricresult1}, we can show more properties of the minimizer $Z=(Z_1, Z_2)$ if $Z_1$ and $Z_2$ are orthogonal on one of the boundaries. The analysis under the polar coordinates plays an important role in the proof.

\begin{theorem}[Theorem \ref{geometryangle} in Section 2]\label{maingeometricresult2}
	Let $Z_1(0), \, Z_1(1)\in \overline{\mathsf{Q}_i} \, $ and $ \, Z_2(0), \, Z_2(1)\in \overline{\mathsf{Q}_j}$, where $\overline{\mathsf{Q}_i}$ and $\overline{\mathsf{Q}_j}$ are two adjacent closed quadrants. Let $Z=(Z_1, \,Z_2) \in H^1([0, \,1], \mathbb{C}^2)$ be a minimizer connecting the two fixed-ends: $Z(0)$ and $Z(1)$. If $Z_1(0) \perp Z_2(0)$ or $Z_1(1)\perp Z_2(1)$, then the polar angles $\theta_1(t)$ and $\theta_2(t)$ of the minimizer $Z=(Z_1, \, Z_2)$ together have at most one critical point for $t\in (0,\,1)$ unless the minimizer is an isosceles triangular motion, a collinear motion or an Euler orbit.
\end{theorem}

\begin{remark}
 By our notation, $Z_1(0) \perp Z_2(0)$ or $Z_1(1)\perp Z_2(1)$ implies that the configuration on one of the fixed-ends could be an isosceles triangle with body 1 as its vertex, an Euler configuration with body 1 at the origin, a binary collision between bodies 2 and 3 or a total collision.
\end{remark}

\begin{remark}
Theorem \ref{maingeometricresult2} implies that at least one of $\theta_i(t) \, (i=1,2)$ is monotone. For example, we can assume $Z_1(t) \in \overline{\mathsf{Q}_4}$,\, $Z_2(t) \in \overline{\mathsf{Q}_3}$ and $\dot{\theta}_2(t) \ne 0$ for all $t \in (0, \,1)$. An illustration is given in Fig. \ref{Jacobifigure3}. It indicates that the region bounded by the curve $Z_2(t)$ and the two half lines $\theta=\theta_2(0)$ and $\theta=\theta_2(1)$ is star-shaped. Furthermore, the curve $Z_1(t)$ can only cross the sector region $C_{[\theta_1(1), \, \theta_1(0)]}$ once.
\end{remark}

In fact, if both fixed ends are isosceles triangles with body 1 as their vertexes, we can show that both polar angles of the minimizer $Z=(Z_1, \, Z_2)$ must be monotone.
\begin{corollary}[Corollary \ref{bothmon} in Section 2]\label{bothmonotone}
Let $Z=(Z_1, \,Z_2) \in H^1([0, \,1], \mathbb{C}^2)$ be a minimizer connecting the two fixed-ends: $Z(0)$ and $Z(1)$. If the polar angles satisfy $\theta_2(0)-\theta_1(0)=\theta_2(1)-\theta_1(1)=\frac{\pi}{2}$, then both $\theta_1$ and $\theta_2$ are monotone. Moreover, the angular momentum of the path $Z=(Z_1, \, Z_2)$ is nonzero.
\end{corollary}

\begin{figure}[htbp]
\begin{center}
\psset{xunit=2.4cm,yunit=2.4cm}
\begin{pspicture}(-2.4,-1)(2.4,0.8)
\psaxes[linewidth=1.5pt]{->}(-1.4, 0.07)(-2.2, -0.9)(-0.5, 0.7)
\psaxes[linewidth=1.5pt]{->}(1.4, 0.07)(0.6, -0.9)(2.3, 0.7)
\psline[linecolor=cyan, doubleline=true]{->}(-0.15, 0.07)(0.15, 0.07)
\rput(-1.3, 0.2){0}
\rput(1.5, 0.2){0}

\psdots[dotsize=2.5pt](-0.8, -0.16)(-1.2, -0.56)(-1.7, -0.66)(-1.9,-0.09)

\rput(-0.5, -0.1){\textcolor{blue}{\small{$Z_1(0)$}}}
\rput(-1, -0.66){\textcolor{blue}{\small{$Z_1(1)$}}}
\rput(-1.9, -0.76){\textcolor{red}{\small{$Z_2(0)$}}}
\rput(-2.1, -0.19){\textcolor{red}{\small{$Z_2(1)$}}}
\psline[linestyle=dashed]{->}(-1.4, 0.07)(-1.9, -0.09)
\psline(-1.525, 0.03)(-1.485,  -0.096 )
\psline(-1.36, -0.056)(-1.485,  -0.096 )
\psline[linestyle=dashed]{->}(-1.4, 0.07)(-1.2, -0.56)
\psline[linewidth=0.4pt, linecolor=green](1.4, 0.07)(0.9, -0.09)( 0.3981,  -0.2506)
\psline[linewidth=0.4pt, linecolor=green](1.4, 0.07)(1.1, -0.66)(1, -0.903)
\psline[linewidth=0.4pt, linecolor=green](1.4, 0.07)(1.6, -0.56)(1.7183,  -0.9327)
\psline[linewidth=0.4pt, linecolor=green](1.4, 0.07)(2, -0.16)(2.3823, -0.3065)

\psdots[dotsize=2.5pt](2, -0.16)(1.6, -0.56)(1.1, -0.66)(0.9,-0.09)
\pscurve[linewidth=2pt, linecolor=blue, linestyle=dotted]{->}(2, -0.16)(1.9,  -0.6)(1.5, -0.8)(1.6, -0.56)
\pscurve[linewidth=2pt, linecolor=red, linestyle=dashed]{->}(1.1, -0.66)(1, -0.56)(0.9,-0.09)
\rput(0.75, -0.5){\textcolor{red}{\small{$Z_2(t)$}}}
\rput(2.1, -0.65){\textcolor{blue}{\small{$Z_1(t)$}}}
 \end{pspicture}
   \end{center}
 \caption{ \label{Jacobifigure3}An illustration of Theorem  \ref{maingeometricresult2}. On the left graph, we assume $Z_1(0), \,Z_1(1) \in \mathsf{\overline{Q}}_4$, $Z_2(0), \, Z_2(1) \in \mathsf{\overline{Q}}_3$ and $Z_1(1) \perp  Z_2(1)$. On the right, $Z_1(t)$ is the dashed blue curve and $Z_2(t)$ is the red one. }
\end{figure}

In the end, we apply Theorem \ref{maingeometricresult1} and Theorem \ref{maingeometricresult2} to two sets of periodic orbits. Theorem \ref{maingeometricresult1} helps exclude possible binary collisions under order constraints in the minimizers and consequently implies their existence, while Theorem \ref{maingeometricresult2} shows interesting geometric properties of the orbits.

The two sets of orbits can be found by considering minimizers connecting collinear configurations and isosceles configurations. Let $Q_S$ be the set of collinear configurations on the $x$-axis with order constraints $q_{2x}(0) \leq q_{1x}(0) \leq q_{3x}(0)$, and $Q_E$ be the set of isosceles triangles, where the symmetry axis of an isosceles is a counterclockwise $\theta$ rotation of the $x$-axis and $q_1$ is  its vertex. Pictures of configurations in $Q_S$ and $Q_E$ are given in Fig. \ref{picQsQe} respectively.

\begin{figure}[htbp]
\begin{center}
\psset{xunit=1.2in,yunit=0.45in}
\begin{pspicture}(-1.5, -0.9)(1.75, 1.3)
\psline[linewidth=0.6pt](-1.5, 0)(-0.15, 0)
\rput(-1.45, 0.25){$q_2$}
\rput(-0.65, 0.25){$q_1$}
\rput(-0.35, 0.25){$q_3$}
\psdot[linecolor=red,dotstyle=*, dotsize=5pt](-1.4,  0)
\psdot[linecolor=blue,dotstyle=*, dotsize=5pt](-0.6, 0)
\psdot[linecolor=black,dotstyle=*, dotsize=5pt](-0.3, 0)
\rput(-1.5, 1.05){$Q_S$:}

\psline[linestyle=dashed, linewidth=0.6pt](0.3, 0)(1.5, 0)
\psline[linestyle=dashed, linewidth=0.6pt](0.7, -0.3)(1.5, 0.5)

\psarc{->}(1, 0){0.6}{0}{22}

\psline[linewidth=0.6pt](0.65,  0.15)(0.78, -0.73)
\psline[linewidth=0.6pt](1.5,  0.5)(0.78, -0.73)
\psline[linewidth=0.6pt](1.5,  0.5)(0.65,  0.15)

\psdot[linecolor=blue,dotstyle=*, dotsize=5pt](1.5,  0.5)
\psdot[linecolor=red,dotstyle=*, dotsize=5pt](0.65,  0.15)
\psdot[linecolor=black,dotstyle=*, dotsize=5pt](0.78, -0.73)

\psline[linecolor=orange]{<-}(1.23, 0.12)(1.68, 0.15)

\rput(0.7, 0.4){$q_2$}
\rput(1.6, 0.55){$q_1$}
\rput(0.88, -0.8){$q_3$}
\rput(1.75, 0.15){$\theta$}

\rput(0.3, 1.05){$Q_E$:}

\end{pspicture}

  \end{center}
\caption{\label{picQsQe} The configurations in $Q_S$ and $Q_E$ are shown, where blue dots represent $q_1$, red dots represent $q_2$ and black dots represent $q_3$. In $Q_S$, three masses are on the $x$-axis with an order $q_{2x} \leq q_{1x} \leq q_{3x}$. In $Q_E$, three masses form an isosceles triangle with $q_1$ as its vertex, while the symmetry axis of the isosceles is a counterclockwise $\theta$ rotation of the $x$-axis. }
\end{figure}

Given $m>0$ and $\theta \in [0, \, \pi/2)$, standard results (For example, Theorem 1.2 in \cite{CH3, FT, Yan3}) imply that there exists a minimizer $\mathcal{P}_{m, \, \theta} \in H^1([0,\,1], \, \Sigma)$, such that
\begin{equation}\label{actionminsetting}
\mathcal{A}(\mathcal{P}_{m, \, \theta})= \inf_{q \in  P(Q_S,\, Q_E)} \mathcal{A}(q) = \inf_{q \in  P(Q_S, \, Q_E)} \int_0^1 (K+U) \, dt,
\end{equation}
where
\[ \mathrm P(Q_S,\, Q_E)=\left  \{q\in H^1([0,\,1],\, \Sigma) \, \big{|}  \, q(0) \in Q_S, \, q(1) \in Q_E \right \}.\]
By the celebrated results of Marchal \cite{Mar} and Chenciner \cite{CA}, $\mathcal{P}_{m, \, \theta}$ is free of collision in $(0,1)$. We are then left to exclude possible boundary collisions. In fact, by following a similar argument in \cite{Yan2, Yan3} and applying Theorem \ref{maingeometricresult1} and Theorem \ref{maingeometricresult2}, we can show that
\begin{proposition}[Proposition \ref{mainthm1} in Section 3]\label{mainthm01}
For each given $m>0$ and $\theta \in [0, \pi/2)$, a minimizer $\mathcal{P}_{m, \, \theta}$ in \eqref{actionminsetting} is collision-free and it can be extended to a periodic or quasi-periodic orbit. Furthermore, if $\mathcal{P}_{m, \, \theta}$ does not coincide with an Euler orbit, then the polar angles $\theta_1(t)$ and $\theta_2(t)$ of the Jacobi coordinates in $\mathcal{P}_{m, \, \theta}$ have at most one critical point.
\end{proposition}

If we change the order of the three masses in the collinear configurations of $Q_S$, it leads to a different set of orbits. Let $Q_{S_1}$ be the set of configurations on the $x$-axis with order constraints $q_{1x}(0) \leq q_{2x}(0) \leq q_{3x}(0)$ (as in Fig. \ref{picQsQe1}). Given $m>0$ and $\theta \in (0, \, \pi/2]$, there exists a minimizer $\widetilde{\mathcal{P}_{m, \, \theta}} \in H^1([0,\,1], \,\Sigma)$, such that
\begin{equation}\label{actionminsetting2}
\mathcal{A}(\widetilde{\mathcal{P}_{m, \, \theta}})= \inf_{q \in  P(Q_{S_1}, \, Q_{E})} \mathcal{A}(q) = \inf_{q \in  P(Q_{S_1}, \, Q_{E})} \int_0^1 (K+U) \, dt,
\end{equation}
where
\[ \mathrm P(Q_{S_1},\, Q_{E})=\left  \{q\in H^1([0,\,1],\, \Sigma) \, \big{|}  \, q(0) \in Q_{S_1}, \, q(1) \in Q_{E} \right \}.\]

\begin{figure}[htbp]
\begin{center}
\psset{xunit=1.2in,yunit=0.45in}
\begin{pspicture}(-1.5, -0.9)(1.75, 1.3)
\psline[linewidth=0.6pt](-1.5, 0)(-0.15, 0)
\rput(-1.45, 0.25){$q_1$}
\rput(-1.05, 0.25){$q_2$}
\rput(-0.35, 0.25){$q_3$}
\psdot[linecolor=blue,dotstyle=*, dotsize=5pt](-1.4,  0)
\psdot[linecolor=red,dotstyle=*, dotsize=5pt](-1, 0)
\psdot[linecolor=black,dotstyle=*, dotsize=5pt](-0.3, 0)
\rput(-1.5, 1.05){$Q_{S_1}$:}

\psline[linestyle=dashed, linewidth=0.6pt](0.4, 0)(1.7, 0)
\psline[linestyle=dashed, linewidth=0.6pt](0.5, -0.5)(1.6, 0.6)

\psarc{->}(1, 0){0.6}{0}{22}

\psline[linewidth=0.6pt](1.64, 0.36)(1.40, 0.65)
\psline[linewidth=0.6pt](0.5, -0.5)(1.64, 0.36)
\psline[linewidth=0.6pt](0.5, -0.5)(1.40, 0.65)

\psdot[linecolor=blue,dotstyle=*, dotsize=5pt](0.5, -0.5)
\psdot[linecolor=red,dotstyle=*, dotsize=5pt](1.40, 0.65)
\psdot[linecolor=black,dotstyle=*, dotsize=5pt](1.64, 0.36)

\psline[linecolor=orange]{<-}(1.23, 0.12)(1.75, 0.15)

\rput(0.62, -0.68){$q_1$}
\rput(1.45, 0.85){$q_2$}
\rput(1.7, 0.56){$q_3$}
\rput(1.8, 0.12){$\theta$}

\rput(0.3, 1.05){$Q_E$:}

\end{pspicture}

  \end{center}
\caption{\label{picQsQe1} The configurations in $Q_{S_1}$ and $Q_E$ are shown. In $Q_{S_1}$, three masses are on the $x$-axis with an order $q_{1x} \leq q_{2x} \leq q_{3x}$. In $Q_E$, three masses form an isosceles triangle with $q_1$ as its vertex, while the symmetry axis of the isosceles is a counterclockwise $\theta$ rotation of the $x$-axis. }
\end{figure}

%
%
%
%
%
%
%
%
%
%
%
%
%
%
%
%
%
%
%

Similarly, it can be shown that
\begin{proposition}[Proposition \ref{retrograderesult} in Section 3]\label{retrograderesult0}
For each given $m>0$ and $\theta \in (0, \, \pi/2)$, a minimizer $\widetilde{\mathcal{P}_{m, \, \theta}}$ in \eqref{actionminsetting2} is collision-free, and it can be extended to a periodic or quasi-periodic orbit. Furthermore, the Jacobi coordinates $Z=(Z_1, \, Z_2)$ of $\widetilde{\mathcal{P}_{m, \, \theta}}$ satisfy that $Z_1 \in \overline{\mathsf{Q}_4}$ and $Z_2 \in \overline{\mathsf{Q}_3}$, which indicates that $\widetilde{\mathcal{P}_{m, \, \theta}}$ contains no collinear configuration except for the boundaries. The corresponding polar angles $\theta_1(t)$ and $\theta_2(t)$ have at most one critical point.
\end{proposition}

\begin{remark}
After extension, $\widetilde{\mathcal{P}_{m, \, \theta}}$ corresponds to a retrograde orbit such that $\phi=4\theta$. Numerically, this retrograde orbit coincide with the one in \cite{CH,CH1}. When $\theta=\pi/2$, it is closely related to one of the open problems \cite{, VEOP} proposed by Venturelli in 2003, in which he asked an existence proof of  the Broucke-H\'{e}non orbit \cite{BR, HM}. By applying Theorem \ref{maingeometricresult1}, the minimizer $\widetilde{\mathcal{P}_{m, \, \pi/2}}$ is shown to be either the Schubart orbit on the $x$-axis, or the Broucke-H\'{e}non orbit in the plane \cite{Yan3}.
\end{remark}

 \begin{figure}[htbp]
    \begin{center}
    \subfigure[ \, $\theta_1$ has a critical point ]{\includegraphics[width=1.6in]{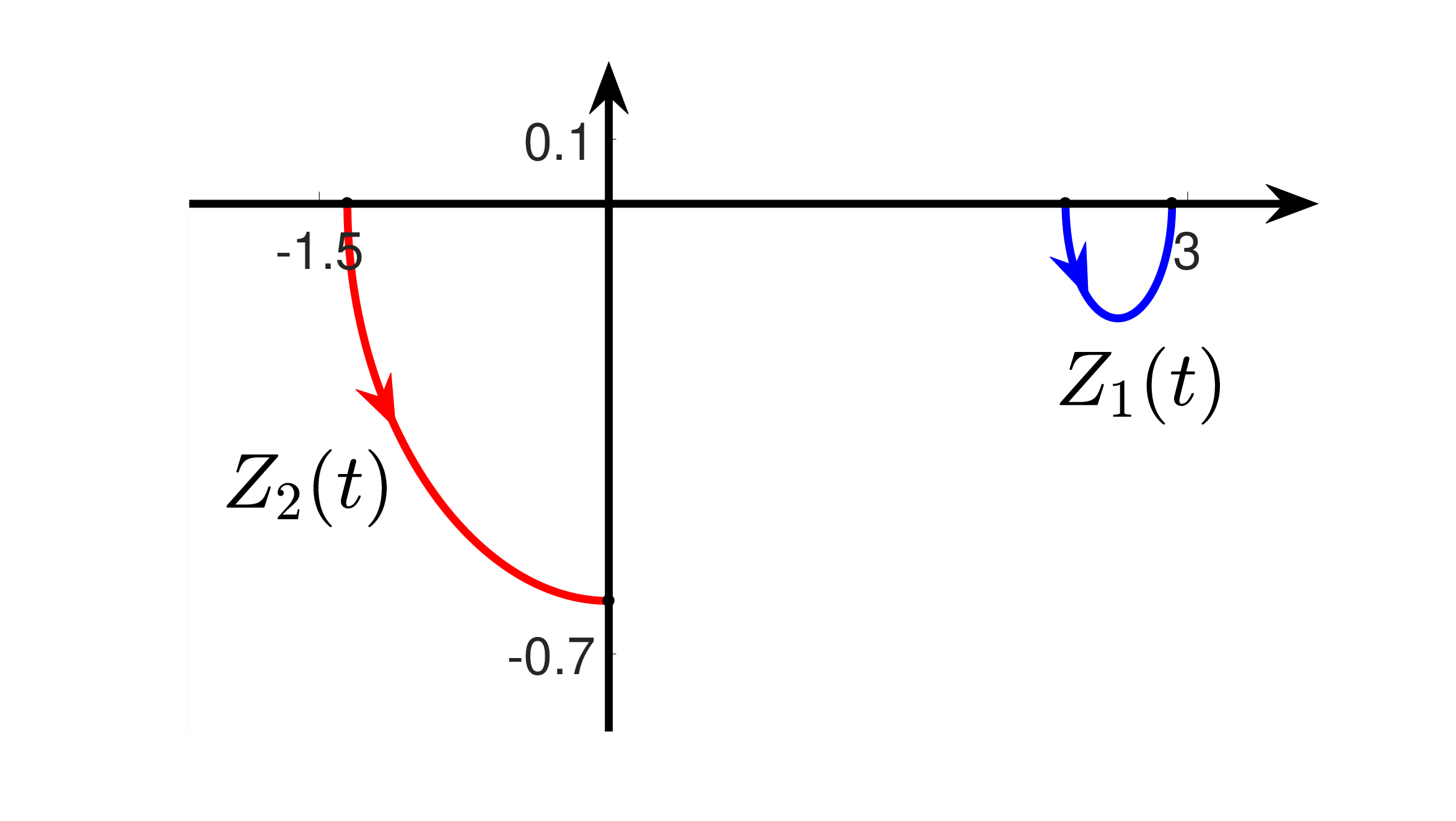}}
    \subfigure[ \, $\theta_2$ has a critical point]{\includegraphics[width=1.6in]{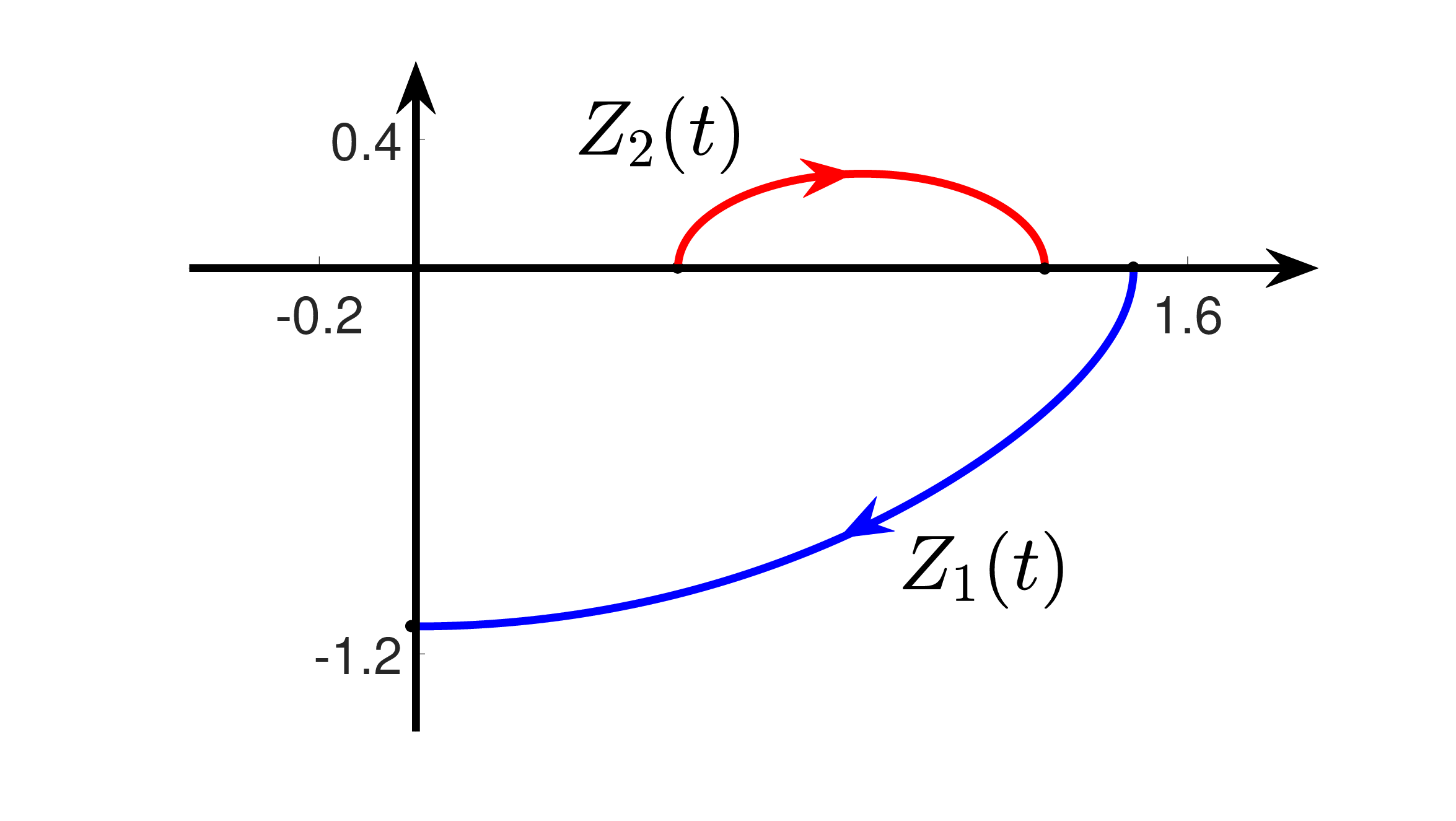}}
    \subfigure[ \, $\theta_1$ and $\theta_2$ have no critical point]{\includegraphics[width=1.9in,height=0.75in]{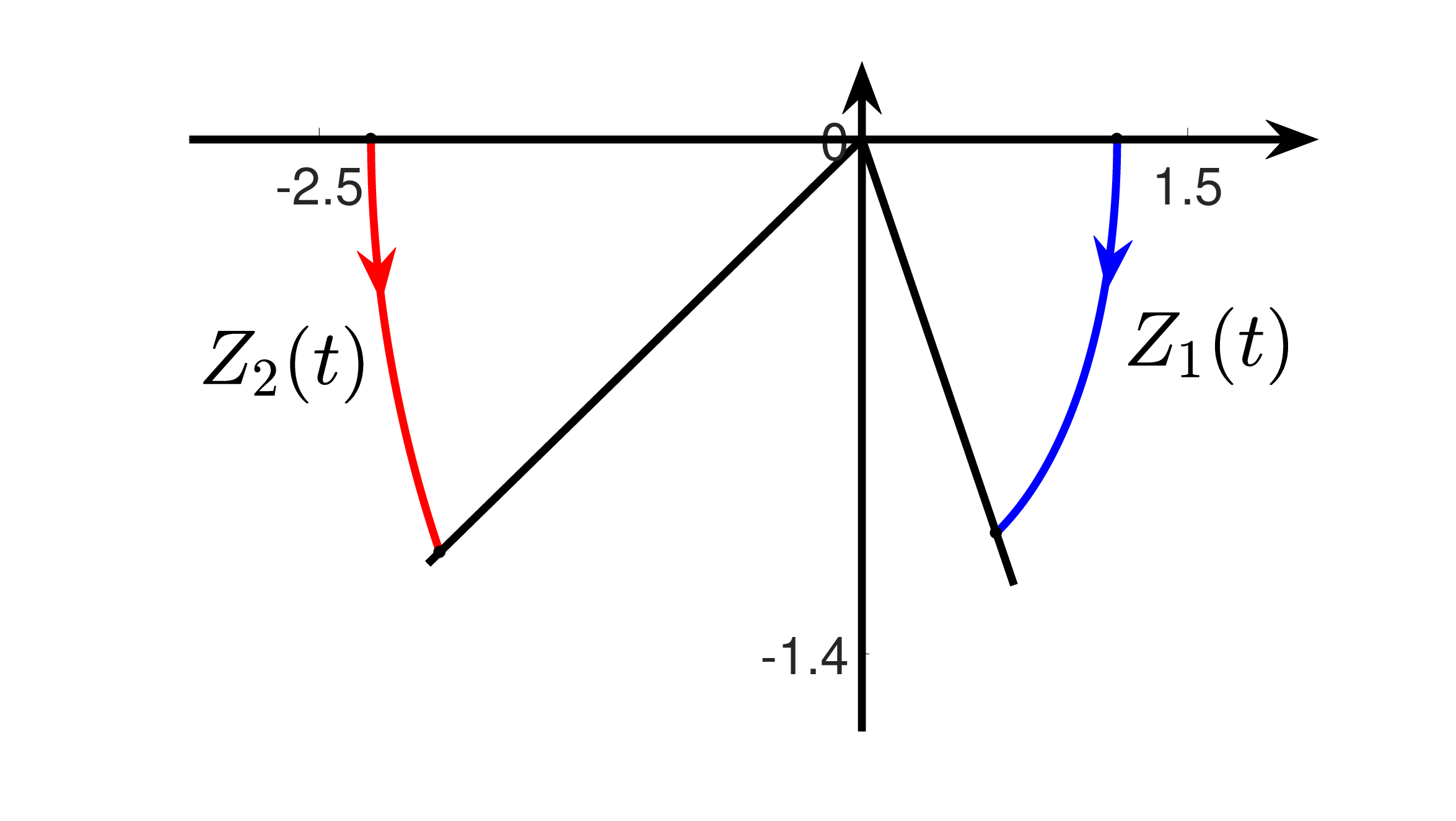}}
   \end{center}
 \caption{ \label{theorem3}The graph of three different minimizers under the Jacobi coordinates, where $Z_1(t)$ is blue and $Z_2(t)$ is red.}
 \end{figure}

\begin{remark}
Theorem \ref{maingeometricresult2} implies three cases for $\theta_i(t) \, (i=1,2)$: $\theta_1(t)$ has a critical point; $\theta_2(t)$ has a critical point; both $\theta_1(t)$ and $\theta_2(t)$ have no critical point. Numerical investigations on the two sets of minimizers $\mathcal{P}_{m, \, \theta}$ and $\widetilde{\mathcal{P}_{m, \, \theta}}$ (as in Fig. \ref{theorem3}) indicate that all three cases can happen. In this sense, Theorem \ref{maingeometricresult2} is a sharp result.
\end{remark}

The paper is organized as follows. In Section \ref{geoarg}, we prove Theorem \ref{maingeometricresult1},  Corollary \ref{corcone}, Theorem \ref{maingeometricresult2} and Corollary \ref{bothmonotone}. In Section \ref{applications}, we apply the results in Section \ref{geoarg} to two sets of periodic orbits in the planar three-body problem.

\section{The Jacobi coordinates and geometric properties of minimizers}\label{geoarg}
In this section, we study properties of minimizers connecting two fixed-ends under the Jacobi coordinates. Recall that $M=[m_1, \, m_2, \, m_3]=[1, \, m, \, m]$ and the Jacobi coordinates are set to be
\[ Z_1=q_3-q_2,  \qquad Z_2=q_1-\frac{q_2+q_3}{2}=(1+\frac{1}{2m})q_1. \]
If we view $Z_1$, $Z_2$ as vectors, we can define $\eta= \arccos \frac{<Z_1, \, Z_2>}{|Z_1||Z_2|}$ to be the angle between $Z_1$ and $Z_2$. It is clear that $\eta \in [0, \pi]$ and the potential $U$ in \eqref{Uformulajacob} can be written in terms of $|Z_1|$, $|Z_2|$ and $\eta$. Note that $U(q), \, U(Z_1,\, Z_2)$ and $U(|Z_1|, \, |Z_2|,\, \eta)$ are the same potential function written in different variables. For simplicity, we used the same notation $U$ to represent them. As in \cite{Yan3}, a direct calculation shows that
\begin{lemma}\label{lemmaU}
The potential function $U=U(|Z_1|,\,|Z_2|,\,\eta)$ is symmetric with respect to $\eta=\frac{\pi}{2}$ and $\eta=0$. $U(|Z_1|,\,|Z_2|,\,\eta)$ is strictly decreasing with respect to $\eta$ when $\eta\in[0,\, \pi/2]$.
\end{lemma}

Let $Z_1=Z_{1x}+ i \, Z_{1y}=r_1 e^{i \theta_1}$ and $Z_2=Z_{2x}+ i \, Z_{2y}=r_2 e^{i \theta_2}$. By \cite{CA, FT, Mar}, the minimizer $Z=(Z_1, \, Z_2)$ is collision-free when $t \in (0, \, 1).$ It follows that it satisfies the equations \eqref{Z1Z2equation0} for $t \in (0, \, 1)$. In the next Lemma, we investigate its motion near $\dot{\theta}_1=0$ or $\dot{\theta}_2=0$.

\begin{lemma}\label{lemmatheta}
	Assume there is some $s \in (0, \,1)$, such that $r_1(s) \ne 0, \, r_2(s) \ne 0$ and $\dot{\theta}_1(s)=0$. Then for $\epsilon>0$ small enough, the motion of $\theta_1$ in $(s-\epsilon, \, s+\epsilon)$ depends on $\theta_2(s)$ and $\dot{\theta}_2(s)$ in the following ways:\\
\hspace*{0.2in}	(a) if $\theta_2(s)\in \big(\theta_1(s)+k\pi, \,  \theta_1(s)+(k+\frac{1}{2})\pi\big)$ for some $k \in \mathrm{Z}$, then $\theta_1(s)$ is a strictly local minimum in $(s-\epsilon, \, s+\epsilon)$;\\
\hspace*{0.2in}	(b) if $\theta_2(s)\in \big(  \theta_1(s)+(k-\frac{1}{2})\pi, \, \theta_1(s)+k\pi \big)$ for some $k \in \mathrm{Z}$, then $\theta_1(s)$ is a strictly local maximum in $(s-\epsilon, \, s+\epsilon)$;\\
\hspace*{0.2in}	(c) if $\theta_2(s)=\theta_1(s)+k\pi$ for some $k \in \mathrm{Z}$, then the motion of $\theta_1$ depends on $\dot{\theta}_2(s)$. When $\dot{\theta}_2(s)=0$, we have $\dot{\theta}_1=\dot{\theta}_2 \equiv 0$ ,which is a collinear motion; when $\dot{\theta}_2(s)>0$, $\theta_1$ is strictly increasing in $(s-\epsilon, \, s+\epsilon)$; when $\dot{\theta}_2(s)<0$, $\theta_1$ is strictly decreasing in $(s-\epsilon, \, s+\epsilon)$; \\
\hspace*{0.2in}	(d) If $\theta_2(s)=\theta_1(s)+(k+\frac{1}{2})\pi$ for some $k \in \mathrm{Z}$, then the motion of $\theta_1$ depends on $\dot{\theta}_2(s)$. When $\dot{\theta}_2(s)=0$, we have $\dot{\theta}_1=\dot{\theta}_2 \equiv 0$ , which is an isosceles triangular motion; when $\dot{\theta}_2(s)>0$, $\theta_1$ is strictly decreasing in $(s-\epsilon, \, s+\epsilon)$; when $\dot{\theta}_2(s)<0$, $\theta_1$ is strictly increasing in $(s-\epsilon, \, s+\epsilon)$.

The same conclusion holds when we exchange $\theta_1$ and $\theta_2$.
\end{lemma}

\begin{proof}
The lemma is shown by analyzing the differential equations in \eqref{Z1Z2equation0}. Note that solutions of \eqref{Z1Z2equation0} are invariant under rotation. Without loss of generality, we can assume $\theta_1(s)=0$. Then $Z_1(s)$ is on the positive $x$-axis.
	
In case (a), $Z_2(t)\in \mathsf{Q}_1\cup \mathsf{Q}_3$. By the first equation in \eqref{Z1Z2equation0},
    \[\ddot{Z}_{1y}(s)>0. \]
Since $\theta_1(s)=\dot{\theta}_1(s)=0$, the polar coordinate form implies that $\ddot{Z}_{1y}(s)=r_1(s)\ddot{\theta}_1(s)$. Thus $\ddot{\theta}_1(s)>0$. We can choose $\epsilon>0$ small enough, such that  $\ddot{\theta}_1(t)>0$ holds for $ t\in (s-\epsilon, \, s+\epsilon)$. Thus  $\theta_1(s)$ is a strictly local minimum in $(s-\epsilon, \,  s+\epsilon)$.
	
In case (b), we have $Z_2(t)\in \mathsf{Q}_2\cup \mathsf{Q}_4$ and $\ddot{\theta}_1(t)<0$ in $(s-\epsilon, \, s+\epsilon)$. Thus $\theta_1(s)$ is a strictly local maximum in $(s-\epsilon,\, s+\epsilon)$.
	
In case (c), \eqref{Z1Z2equation0} implies that $\ddot{Z}_{1y}(s)=0$. Note that $\ddot{Z}_{1y}(s)=r_1(s)\ddot{\theta}_1(s)$, it implies that $\theta_1(s)=\dot{\theta}_1(s)=\ddot{\theta}_1(s)=0$ and $Z_{1y}(s)=Z_{2y}(s)=\dot{Z}_{1y}(s)=0$. We then consider the position $Z_2(t)$ around $t=s$. If $\dot{\theta}_2(s)=0$, by the existence and uniqueness theorem of ODE equations, it must be a collinear motion. If $\dot{\theta}_2(s)>0$, we consider the sign of $\dddot{\theta}_1(s)$. In fact, $\dddot{Z}_{1y}(s)=r_1(s) \dddot{\theta}_1(s)$. On the other hand, by differentiating the first equation in \eqref{Z1Z2equation0}, it follows that
\begin{equation}\label{dddz1y}
\dddot{Z}_{1y}(s)=   \frac{\dot{Z}_{2y}(s)}{|Z_2(s)-Z_1(s)/2|^3}- \frac{\dot{Z}_{2y}(s)}{|Z_2(s)+Z_1(s)/2|^3}.
\end{equation}
 Note that $\theta_2=k \pi$ and $\dot{Z}_{2y}(s)= r_2(s) \cos (k \pi) \dot{\theta}_2(s)$. By \eqref{dddz1y}, $\dddot{Z}_{1y}(s)=r_1(s) \dddot{\theta}_1(s)$ has the same sign as $\dot{\theta}_2(s)$. Hence, there exists $\epsilon>0$ small enough, such that $\theta_1$ is strictly increasing in $(s-\epsilon, \, s+\epsilon)$ when $\dot{\theta}_2(s)>0$. Similarly, $\theta_1$ is strictly decreasing in $(s-\epsilon, \, s+\epsilon)$ when $\dot{\theta}_2(s)<0$.

In case (d), note that $\ddot{Z}_{1y}(s)=r_1(s)\ddot{\theta}_1(s)$, $\theta_1(s)=\dot{\theta}_1(s)=\ddot{\theta}_1(s)=0$ and $Z_{1y}(s)=Z_{2x}(s)=\dot{Z}_{1y}(s)=0$. By differentiating the first equation in \eqref{Z1Z2equation0}, it follows that
\begin{equation}\label{dddz1y1}
\dddot{Z}_{1y}(s)=   \frac{3Z_{2y}(s) Z_{1x}(s) \dot{Z}_{2x}(s)}{|Z_2(s)-Z_1(s)/2|^5}.
\end{equation}
Since $Z_{2y}(s) Z_{1x}(s) \dot{Z}_{2x}(s)= -r_1(s) r_2^2(s) \dot{\theta}_2(s)$, it implies that $\dddot{\theta}_1(s)$ and $\dot{\theta}_2(s)$ have opposite signs. Hence, case (d) holds.

When we exchange $\theta_1$ and $\theta_2$, the proof follows by a similar argument.
\end{proof}

By Lemma \ref{lemmaU} and Lemma \ref{lemmatheta}, we can show that
\begin{theorem}[Theorem \ref{maingeometricresult1}]\label{thmquadrants}
Assume $Z_1(0), \, Z_1(1)\in \overline{\mathsf{Q}_i} \, $ and $ \, Z_2(0), \, Z_2(1)\in \overline{\mathsf{Q}_j}$, where $\overline{\mathsf{Q}_i}$ and $\overline{\mathsf{Q}_j}$ are two adjacent closed quadrants. Let $Z=(Z_1, \,Z_2) \in H^1([0, \,1], \mathbb{C}^2)$ be a minimizer connecting the two fixed-ends: $Z(0)=(Z_1(0), \, Z_2(0))$ and $Z(1)=(Z_1(1), \, Z_2(1))$. That is,
\begin{equation}\label{Bolza2}
 \mathcal{A}(Z)=\inf_{\widehat{Z} \in P(Z(0), \, Z(1))}\mathcal{A}(\widehat{Z})= \inf_{\widehat{Z} \in P(Z(0), \, Z(1))} \int_0^1 (K+U) \, dt,
 \end{equation}
where $K$ (in \eqref{Kformulajacob}) and $U$ (in \eqref{Uformulajacob}) are the kinetic energy and potential function respectively, and
\[P(Z(0), \, Z(1))= \left\{\, \widehat{Z}=(\widehat{Z}_1, \, \widehat{Z}_2) \in H^1([0,1], \mathbb{C}^2) \, \big| \, \widehat{Z}(0)=Z(0), \, \,  \widehat{Z}(1)=Z(1) \, \right\}. \]
Then $Z_1(t)$ and $Z_2(t)$ must be in two adjacent closed quadrants for all $t \in [0,\, 1]$ and they must satisfy one of the following three cases:\\
\hspace*{0.2in} (a) both $Z_1(t)$ and $Z_2(t)$ are away from the coordinate axes for all $t \in (0, \, 1)$; \\
\hspace*{0.2in}  (b) both $Z_1(t)$ and $Z_2(t)$ are on the coordinate axes for all $t\in [0, \,1]$ and they never touch the origin in $(0,\, 1)$;  \\
\hspace*{0.2in}  (c) $Z=(Z_1, \, Z_2)$ is part of an Euler orbit with $Z_2(t) \equiv 0$ for all $t \in [0,\,  1]$.
\end{theorem}

\begin{proof}
	 Without loss of generality, we assume $Z_1(0), \, Z_1(1)\in \overline{Q}_1$ and $Z_2(0), \, Z_2(1)\in \overline{Q}_2$. Let $\widetilde{Z}=(\widetilde{Z_1},\, \widetilde{Z_2})$ be
	 \[\widetilde{Z_1}(t)=\big(|Z_{1x}(t)|, \, |Z_{1y}(t)|\big), \qquad  \widetilde{Z_2}(t)=\big(-|Z_{2x}(t)|, \, |Z_{2y}(t)|\big).\]
By assumption, the two paths $Z=(Z_1,\,Z_2)$ and $\widetilde{Z}=(\widetilde{Z_1},\, \widetilde{Z_2})$ have the same boundaries. For $\forall t\in[0,\,1]$, $\widetilde{Z_1}(t)\in \overline{Q}_1$ and $\widetilde{Z_2}(t)\in \overline{Q}_2$.

Note that $|Z_1(t)|=|\widetilde{Z_1}(t)|$ and $|Z_2(t)|=|\widetilde{Z_2}(t)|$. By comparing the angles $\eta(t)$ and $\tilde{\eta}(t)$ of the two paths, Lemma \ref{lemmaU} implies that
	 \begin{equation}\label{UZ}
	 U(Z(t))\ge U(\widetilde{Z}(t)).
	 \end{equation}
The equality holds if and only if $Z_1(t)$ and $Z_2(t)$ are in two adjacent closed quadrants.
	
By the definition of $\tilde{Z}$, it follows that the integral of the kinetic energy $\displaystyle \int_0^1 K \, dt$ is the same for the two paths: $Z$ and $\tilde{Z}$. Hence,
	 \begin{equation}\label{AZ}
	\mathcal{A}(Z) \ge \mathcal{A}(\widetilde{Z}).
	 \end{equation}
Since $Z=(Z_1, \, Z_2)$  is a minimizer connecting two fixed-ends and $\widetilde{Z}=(\widetilde{Z_1}, \, \widetilde{Z_2})$ is an admissible path with the same boundaries, the equality in \eqref{AZ} must holds. Both $Z$ and $\widetilde{Z}$ are minimizers of the same fixed-ends problem.
	
	 By [11] and [2], the minimizers $Z$ and $\widetilde{Z}$ are collision-free in $(0, \,1)$, and they correspond to solutions of the Newtonian equations \eqref{newton} in $(0,\, 1)$. Both $Z(t)$ and $\widetilde{Z}(t)$ must be analytic for $ t \in (0,\,1)$. Thus $Z_1$ and $Z_2$ can not cross the coordinate axes. In other words, $Z_1$ and $Z_2$ are always in two adjacent closed quadrants. If $Z_1$ and $Z_2$ do not touch the coordinate axes in $(0,\,1)$, it's case $(a)$ in the theorem. If $Z_1$ or $Z_2$ touches the coordinate axes in $(0,\,1)$, it must be tangent to the coordinate axes. We then prove that $Z=(Z_1, \, Z_2)$ must be either case $(b)$ or case $(c)$ in the theorem.
	
	Note that $Z_1=q_3-q_2$ and the solution is collision-free in $(0,\,1)$, it follows that $Z_1(t) \ne 0$ for all $t\in(0,\,1)$. If there exists some $t_0\in (0,\,1)$ such that $Z_2(t_0)=0$, by the analyticity of $Z_2(t)$ and $\widetilde{Z_2}(t)$, we have $\dot{Z}_2(t_0) = 0$. Since $\{(q_1, \, q_2, \, q_3) \in \Sigma\, | \, q_1=0,  \, \dot{q}_1= 0\}$ is an invariant set, it implies that $Z_2(t)\equiv 0$ for all $t\in [0,\, 1]$. In this case, the motion is part of an Euler orbit, which is case $(c)$.
	
	Assume that $Z_1(t_0)$ or $Z_2(t_0)$ is on the axes away from the origin for some $t_0 \in (0,\, 1)$. It follows that $\dot{\theta}_1(t_0)=0$ or $\dot{\theta}_2(t_0)=0$. Without loss of generality, we assume $\dot{\theta}_1(t_0)=0$. If cases $(a)$ or $(b)$ of Lemma \ref{lemmatheta} happens, it contradicts the fact that $Z_1(t)$ and $Z_2(t)$ are always in two adjacent closed quadrants. Hence $Z_2(t_0)$ is on the axes. If $\dot{\theta}_2(t_0) \ne 0$, by cases $(c)$ and $(d)$ of Lemma \ref{lemmatheta}, $Z_1$ will cross the axes. Contradiction! Thus there must be $\dot{\theta}_2(t_0)=0$. By Lemma \ref{lemmatheta}, it is a collinear motion or an isosceles triangular motion, which is case $(b)$.	

\end{proof}

By considering the minimizer $Z=(Z_1, \, Z_2)$ in two different coordinate systems, we can extend Theorem \ref{thmquadrants} to the following Corollary.
\begin{corollary}[Corollary \ref{corcone}]\label{thmcone}
		Assume there exist $a, b \in \mathbb{R}$ satisfying $0\le b-a<\frac{\pi}{2}$, such that $Z_1(0)$, \,$Z_1(1)$ $\in C_{[a,\,b]}$, $Z_2(0)$,\, $Z_2(1)\in C_{[a,\,b]}^+$. We further assume that $|Z_i(0)|+|Z_i(1)| \ne 0, \, (i=1,\,2)$. Let $Z=(Z_1, \,Z_2) \in H^1([0, \,1], \mathbb{C}^2)$ be a minimizer connecting the two fixed-ends: $Z(0)$ and $Z(1)$.

Then $Z_1(t)\in C_{[a,\,b]}$ and $Z_2(t)\in C_{[a,\,b]}^+$ for all $t\in [0,\,1]$, and they must satisfy one of the following two cases:\\
\hspace*{0.4in}	    (a) both $Z_1(t)$ and $Z_2(t)$ are away from the boundaries of the corresponding closed cones for all $ t \in (0, \,1)$;\\
\hspace*{0.4in} 	(b) both $Z_1(t)$ and $Z_2(t)$ are on the boundaries of the corresponding closed cones for all $t\in [0, \,1]$, and they can not reach the origin when $ t \in (0, \,1)$.

Similar conclusion holds if $Z_2(0), \, Z_2(1)\in C_{[a,\,b]}^-$.
\end{corollary}

\begin{proof}
We only prove the case when $Z_1(0), \,Z_1(1)\in C_{[a,\,b]}$ and $Z_2(0), \,Z_2(1)\in C_{[a,\,b]}^+$, while the other case when $Z_1(0), \,Z_1(1)\in C_{[a,\,b]}$ and $Z_2(0), \, Z_2(1)\in C_{[a,\,b]}^-$ follows similarly. According to the boundary condition, we can consider the following two coordinate systems:\\
\hspace*{0.3in} (i) the $x$-axis coincides with the line $\theta=a$ and the $y$-axis coincides with the line $\theta=a+\frac{\pi}{2}$; \\
\hspace*{0.3in} (ii) the $x$-axis coincides with the line $\theta=b$ and the $y$-axis coincides with the line $\theta=b+\frac{\pi}{2}$.

Note that the two coordinate systems coincide when $a=b$.

 By assumption, case $(c)$ of Theorem \ref{thmquadrants} can not occur. If case $(b)$ of Theorem \ref{thmquadrants} happens, it must be an isosceles triangular motion and it satisfies that $Z_1(t)\in C_{[a,\,b]}$ and $Z_2(t)\in C_{[a,\,b]}^+$.

 If case $(a)$ of Theorem \ref{thmquadrants} happens, $Z_1(t)$ is away from the coordinate axes for all $t\in (0,\,1)$. We shall prove that $\theta_1(t)\in [a, \, b]$ for $\forall t\in [0,\, 1]$. Note that the coordinate axes of the two coordinate systems divide the plane into eight closed cones defined by eight intervals of angles (when $a=b$, it's four cones). By assumption, we may assume $Z_1(0) \ne 0$ and $Z_2(0) \ne 0$. Since $Z_1(t)$ can not touch the coordinate axes in $(0,1)$, there are three situations: (i) $\theta_1(t)\in [a,\, b]$ for $\forall t\in [0,\, 1)$; (ii) $\theta_1(t)\in [b, \, a+\frac{\pi}{2}]$ for $\forall t\in [0, \, 1)$; (iii) $\theta_1(t) \in [b-\frac{\pi}{2}, \, a]$ for $\forall t\in [0,\, 1)$. (For simplicity, we still use the same notation $b-\frac{\pi}{2}$ for $b<\frac{\pi}{2}$.)
	
We first show that (ii) can not happen in the case when $Z_1(1)$ and $Z_2(1)$ are nonzero. Assume $\theta_1(t)\in [b, \, a+\frac{\pi}{2}]$, then $\theta_1(0)=\theta_1(1)=b$. By applying Theorem \ref{thmquadrants} in coordinate system (ii), it follows that $\theta_1(t)\in [b,\, b+\frac{\pi}{2}]$ and $\theta_2(t)\in [b+\frac{\pi}{2},\, b+\pi]$. Thus $\theta_2(0)=\theta_2(1)=b+\frac{\pi}{2}$. For this case, we claim that the minimizer satisfies $\dot{\theta}_1=\dot{\theta}_2= 0$. If not, we define a new path $\widetilde{Z}=(\widetilde{Z}_1, \, \widetilde{Z}_2)$: $\widetilde{Z}_1(t)=r_1(t) e^{i b}$ and $\widetilde{Z}_2(t)=r_2(t) e^{i (b+\pi/2)}$. It is clear that $K(\widetilde{Z})<K(Z)$. By Lemma \ref{lemmaU}, $U(\widetilde{Z}) \leq U(Z)$. It follows that $\mathcal{A}(\widetilde{Z})<\mathcal{A}(Z)$. It contradicts the fact that $Z$ is a minimizer connecting the two fixed-ends. Hence, the motion is an isosceles triangular motion with $\dot{\theta}_1=\dot{\theta}_2 \equiv 0$, which is case $(b)$ of Theorem \ref{thmquadrants}. Contradiction to the assumption that case $(a)$ of Theorem \ref{thmquadrants} happens! Indeed, if $Z_1(1)=0$ or $Z_2(1)=0$, one can find a contradiction by the same argument. Similarly, the case when $\theta_1(t)\in [b-\frac{\pi}{2},\,a]$ can be excluded.

 Thus $\theta_1(t)\in [a,\,b]$ for $\forall t\in [0,\,1)$. That is $Z_1(t)\in C_{[a,\,b]}$ for any $t\in [0,\,1)$.
\end{proof}

Now we can apply Lemma \ref{lemmatheta} and Corollary \ref{thmcone} to prove the following theorem.

\begin{theorem}[Theorem \ref{maingeometricresult2}]\label{geometryangle}
	Let $Z_1(0), \, Z_1(1)\in \overline{\mathsf{Q}_i} \, $ and $ \, Z_2(0), \, Z_2(1)\in \overline{\mathsf{Q}_j}$, where $\overline{\mathsf{Q}_i}$ and $\overline{\mathsf{Q}_j}$ are two adjacent closed quadrants. Let $Z=(Z_1, \,Z_2) \in H^1([0, \,1], \mathbb{C}^2)$ be a minimizer connecting the two fixed-ends: $Z(0)$ and $Z(1)$. If $Z_1(0) \perp Z_2(0)$ or $Z_1(1)\perp Z_2(1)$, then the polar angles $\theta_1(t)$ and $\theta_2(t)$ of the minimizer $Z=(Z_1, \, Z_2)$ together have at most one critical point for $t\in (0,\,1)$ unless the minimizer is an isosceles triangular motion, a collinear motion or an Euler orbit.
\end{theorem}
\begin{proof}
	Without loss of generality, we assume $Z_1 \in \overline{\mathsf{Q}_4}$, \, $Z_2 \in \overline{\mathsf{Q}_3}$ and $Z_1(1) \perp Z_2(1)$. We should further assume that the minimizer $Z=(Z_1, \, Z_2)$ is not an isosceles triangular motion, a collinear motion or an Euler orbit. By Theorem \ref{thmquadrants}, $Z_1(t)\in Q_4, \, Z_2(t)\in Q_3$ for $t\in (0,\,1)$. Hence, both $r_1(t)$ and $r_2(t)$ are nonzero for $t \in (0,\,1)$ and Lemma \ref{lemmatheta} can be applied. The idea is to analyze the cases in Lemma \ref{lemmatheta} and derive contradictions with Corollary \ref{thmcone} if the polar angles have more than one critical point. We show this theorem in three steps.
	
	\textbf{Step 1:} Prove that cases $(c)$ and $(d)$ of Lemma \ref{lemmatheta} can not happen for $t \in (0,1)$.
	
	Since $Z_1 \in \overline{\mathsf{Q}_4}$, \, $Z_2 \in \overline{\mathsf{Q}_3}$, it follows that the only possible case when $Z_1(t)$ and $Z_2(t)$ are collinear is when they are on the $x$-axis. If case $(c)$ of Lemma \ref{lemmatheta} happens at $t_0\in (0,1)$, Theorem \ref{maingeometricresult1} implies that $Z_1$ and $Z_2$ stay on the $x$-axis for all $t \in [0,1]$. Contradiction to the assumption!
	
   If case $(d)$ happens, then $\theta_2(t_0)+\frac{\pi}{2}=\theta_1(t_0)$ for some $t_0 \in (0, \, 1)$. The following three cases will be considered.

\textsl{Case (i):} If both $Z_1(1)$ and $Z_2(1)$ are nonzero, then the two configurations at $t=t_0$ and $t=1$ satisfy the assumption in Corollary \ref{thmcone}. It follows that $Z_1(t) \in C_{[\theta_1(t_0), \, \theta_1(1)]}$ and $Z_2(t) \in C_{[\theta_2(t_0), \, \theta_2(1)]}$ for $t\in [t_0,\,1]$. (The two cones could be $C_{[\theta_1(1), \, \theta_1(t_0)]}$ and $C_{[\theta_2(1), \, \theta_2(t_0)]}$ if $\theta_1(1)<\theta_1(t_0)$.) On the other hand,  since the motion is not an isosceles triangular motion, Lemma \ref{lemmatheta} implies that $\dot{\theta}_2(t_0) \ne 0$, and $\dot{\theta}_1$ and $\dot{\theta}_2$ have opposite monotonicity in $(t_0,\, t_0+\epsilon)$ for some $\epsilon>0$ small enough. Contradiction!

 \textsl{Case (ii):}  If $Z_1(1)=0$ and $Z_2(1) \ne 0$, Corollary \ref{thmcone} implies that $Z_1(t) \in C_{[\theta_1(t_0), \, \theta_2(1)+ \pi/2]}$ and $Z_2(t) \in C_{[\theta_2(t_0),\, \theta_2(1)]}$ for $t\in [t_0,\,1]$. (The two cones could be $C_{[\theta_2(1)+ \pi/2, \, \theta_1(t_0)]}$ and $C_{[\theta_2(1),\, \theta_2(t_0)]}$ if $\theta_2(1)<\theta_2(t_0)$.) However, $\dot{\theta}_1$ and $\dot{\theta}_2$ have opposite monotonicity in $(t_0,\, t_0+\epsilon)$ for some $\epsilon>0$ small enough. Contradiction! Similarly, a contradiction can be found when $Z_2(1)=0$ and $Z_1(1) \ne 0$.

\textsl{Case (iii):}   If $Z_1(1)=Z_2(1)= 0$, Corollary \ref{thmcone} implies that $\dot{\theta}_1(t)=\dot{\theta}_2(t)=0$ for $t \in [t_0, \, 1]$. The minimizer then has an isosceles triangular motion. Contradiction!

    Hence, cases $(c)$ and $(d)$ of Lemma \ref{lemmatheta} can not happen.
	
	\textbf{Step 2:} Prove that $\theta_2(t)$ has at most one critical point for $t \in (0, \, 1)$. So does $\theta_1(t)$.
	
	Assume $\theta_2$ has at least two critical point. Let $0<a<b<1$ be the two largest points such that $\dot{\theta}_2(a)=\dot{\theta}_2(b)=0$. By \textbf{Step 1} and Lemma \ref{lemmatheta}, one point is a local minimum of $\theta_2$ and the other point is a local maximum. Without loss of generality, we assume $\theta_2(a)$ is a local minimum and $\theta_2(b)$ is a local maximum. Note that $Z_1 \in \overline{\mathsf{Q}_4}$, $Z_2 \in \overline{\mathsf{Q}_3}$, it implies that $\theta_1-\theta_2 \in [0,\, \pi]$. By Lemma \ref{lemmatheta},
	\begin{equation}\label{thetaab}
	\begin{split}
	0<\theta_1(a)-\theta_2(a)<\frac{\pi}{2},\\
	\frac{\pi}{2}<\theta_1(b)-\theta_2(b)<\pi.\\
	\end{split}
	\end{equation}
	By an intermediate value theorem, there exists some $s\in (a, \, b)$ such that $\theta_1(s)-\theta_2(s)=\frac{\pi}{2}$. There are three cases to be considered.
	
	\textsl{Case (i):} If $Z_2(1) \ne 0$ and $Z_1(1) \ne 0$, then $Z_1(1) \perp Z_2(1)$ implies that $\theta_1(1)-\theta_2(1)=\frac{\pi}{2}$. By Corollary \ref{thmcone}, it follows that $Z_2(t) \in C_{[\theta_2(s),\theta_2(1)]}$ holds for $t \in [s, \, 1]$. Since the local maximum $\theta_2(b)$ is the only critical point of $\theta_2$ in $(s,1)$, it implies that
	\begin{equation}
	\theta_2(b)>\theta_2(s), \qquad  \theta_2(b)>\theta_2(1).
	\end{equation}
	Contradiction to the fact that $Z_2(t)\in C_{[\theta_2(s), \, \theta_2(1)]}$ for $t\in [s,1]$!
	
	\textsl{Case (ii):} If $Z_2(1)=0$ and $Z_1(1) \ne 0$, then $\theta_2(1)$ is not well defined. Consider the motion in the time interval $t \in [a,\, 1]$. Note that $0<\theta_1(a)-\theta_2(a)<\frac{\pi}{2}$. If $\theta_1(1)-\theta_2(a)<\frac{\pi}{2}$, then Corollary \ref{thmcone} implies that the minimizing path satisfies $Z_1(t)\in C_{[\theta_1(a), \, \theta_2(a)+\frac{\pi}{2}]}$ and $Z_2(t)\in C_{[\theta_1(a)-\frac{\pi}{2}, \, \theta_2(a)]}$ for $t \in [a,\, 1]$. This contradicts the fact that $\theta_2(a)$ is a local minimum. Thus there must have $\theta_1(1)-\theta_2(a)\ge\frac{\pi}{2}$.
	Similarly, by considering the interval $t \in [b,\,1]$, we have $\theta_1(1)-\theta_2(b)\le \frac{\pi}{2}$. The above argument shows that
	\begin{equation}
	\theta_2(a)+\frac{\pi}{2}\le\theta_1(1)\le \theta_2(b)+ \frac{\pi}{2}.
	\end{equation}
	Now we have $Z_1(s),\, Z_1(1)\in C_{[\theta_2(a)+\frac{\pi}{2}, \, \theta_2(b)+\frac{\pi}{2}]}$ and $Z_2(s),\,Z_2(1)\in C_{[\theta_2(a),\, \theta_2(b)]}$. By Corollary \ref{thmcone}, $Z_1(t)\in C_{[\theta_2(a)+\frac{\pi}{2},\, \theta_2(b)+\frac{\pi}{2}]}$ for any $t\in [s,\,1]$. It follows that $\theta_1(b)\le \theta_2(b)+ \pi/2$. Contradiction to \eqref{thetaab}!

	If $Z_1(1)=0$ and $Z_2(1)\ne 0$, a contradiction can be found by a similar argument.

\textsl{Case (iii):} If $Z_2(1)=Z_1(1)=0$, then neither $\theta_2(1)$ nor $\theta_1(1)$ is well defined. By applying Corollary \ref{thmcone}, it follows that for $\forall t\in [s,\,1]$, $\dot{\theta}_1(t)=\dot{\theta}_2(t)=0$. The minimizer has an isosceles triangular motion. Contradiction!
	
Therefore, there is at most one point $t_0 \in (0,\,1)$, such that $\dot{\theta}_2( t_0)=0$. Similarly, one can show that $\theta_1(t)$ has at most one critical point.
	
	\textbf{Step 3:} It's impossible that both $\theta_1(t)$ and $\theta_2(t)$ have critical points in $(0,\,1)$.
	
	We show it by contradiction. Assume that $\dot{\theta}_1(a')=\dot{\theta}_2(b')=0$. Without loss of generality, we can assume $0<a' \le b'<1$. By \textbf{Step 1}, they must be either local maximum or local minimum. By \textbf{Step 2}, they are the only critical point for $\theta_1$ and $\theta_2$ in $(0, \, 1)$.

	If both $\theta_1(a')$ and $\theta_2(b')$ are local maximum, then $\theta_1$ is decreasing in $[a',1)$, $\theta_2$ is decreasing in $[b',1)$ and increasing in $[a',b']$. Since $Z_1(t)\in Q_4$ and $Z_2(t)\in Q_3$, Lemma \ref{lemmatheta} implies that
	\begin{equation}
	0<\theta_1(a')-\theta_2(a')<\frac{\pi}{2},\qquad  \theta_1(b')-\theta_2(b')>\frac{\pi}{2}.
	\end{equation}
	This is impossible because $\theta_1(a')-\theta_2(a') \ge \theta_1(b')-\theta_2(b')$. Similarly, it's also impossible that both $\theta_1(a')$ and $\theta_2(b')$ are local minimum.	
	
	If $\theta_1(a')$ is a local minimum and $\theta_2(b')$ is a local maximum, then $\theta_1$ is increasing in $[a',1)$, $\theta_2$ is decreasing in $[b',1)$ and increasing in $[a',b']$. By Lemma \ref{lemmatheta},
	\begin{equation}\label{step3.2}
	\theta_1(a')-\theta_2(a')>\frac{\pi}{2}, \qquad \theta_1(b')-\theta_2(b')>\frac{\pi}{2}.
	\end{equation}
	In this situation, the proof is divided into the following three cases.
	
	\textsl{Case (i):} If $Z_2(1)\ne 0$ and $Z_1(1)\ne 0$, the assumption $Z_1(1)\perp Z_2(1) $ implies that $\theta_1(1)-\theta_2(1)=\frac{\pi}{2}$. It follows that $\theta_1(b')-\theta_2(b')<\theta_1(1)-\theta_2(1)= \pi/2$. Contradiction to \eqref{step3.2}!
	
	\textsl{Case (ii):} If $Z_2(1)\ne 0$ and $Z_1(1)=0$, then $\theta_2(1)$ is well defined. Define $\theta_1(1)=\theta_2(1)+\frac{\pi}{2}$.
	
	When $\theta_1(a') \le \theta_1(1)$, then $\theta_1(1) \ge \theta_1(a') > \theta_2(a')+\pi/2$ by \eqref{step3.2}. Corollary \ref{thmcone} implies that for $\forall t\in [a',1]$,
	\begin{equation}
	Z_1(t)\in C_{[\theta_2(a')+\frac{\pi}{2},\, \theta_1(1)]}, \ \ \ Z_2(t)\in C_{[\theta_2(a'),\, \theta_2(1)]} .
	\end{equation}
	This contradicts the assumption that $\theta_2(b')>\theta_2(1)$.
	
	When $\theta_1(a')> \theta_1(1)$, we consider the motion when $t \in [a', \, 1]$. Note that $Z_1(t)\in Q_4$ and $Z_2(t)\in Q_3$. By Corollary \ref{thmcone},
\[ 	Z_1(t)\in C_{[\frac{3 \pi}{2},  \, \theta_1(a')]}, \qquad Z_2(t)\in C_{[\pi, \, \theta_1(a')-\frac{\pi}{2}]},  \quad \text{for} \, \, \, \forall t \in [a', \, 1].  \]
It follows that
	\begin{equation}
	\theta_1(t)\le \theta_1(a'), \qquad \text{for} \, \, \,  \forall t\in [a',1].
	\end{equation}
	This contradicts the fact that $\theta_1(a')$ is a local minimum.
	
	When $Z_2(1)=0$ and $Z_1(1)\ne 0$, a contradiction can be found by a similar argument.
	
	\textsl{Case (iii):} If $Z_2(1)=Z_1(1)=0$, then neither $\theta_2(1)$ nor $\theta_1(1)$ is well defined. By applying Corollary \ref{thmcone} to the time interval $[a',1]$, we have for $\forall t\in [a',1]$,
	\begin{equation}
	Z_1(t)\in C_{[\theta_2(a')+\frac{\pi}{2}, \, \theta_1(a')]}, \ \ \ Z_2(t)\in C_{[\theta_2(a'),\, \theta_1(a')-\frac{\pi}{2}]} .
	\end{equation}
	This contradicts the fact that $\theta_1(a')$ is a local minimum.
	
Hence, the case when $\theta_1(a')$ is a local minimum and $\theta_2(b')$ is a local maximum can not happen. For the case when $\theta_1(a')$ is a local maximum and $\theta_2(b')$ is a local minimum, contradictions can be found by a similar argument.

Therefore, $\theta_1(t)$ and $\theta_2(t)$ together have at most one critical point in $(0,\,1)$.
\end{proof}

By applying Corollary \ref{thmcone} and Theorem \ref{geometryangle}, we can prove the following corollary.
\begin{corollary}[Corollary \ref{bothmonotone}]\label{bothmon}
Let $Z=(Z_1, \,Z_2) \in H^1([0, \,1], \mathbb{C}^2)$ be a minimizer connecting the two fixed-ends: $Z(0)$ and $Z(1)$. If the polar angles satisfy $\theta_2(0)-\theta_1(0)=\theta_2(1)-\theta_1(1)=\frac{\pi}{2}$, then both $\theta_1$ and $\theta_2$ are monotone. Moreover, the angular momentum of the path $Z=(Z_1, \, Z_2)$ is nonzero. 
\end{corollary}
\begin{proof}
Let $\theta_1(0)=a$, $\theta_1(1)=b$. Without loss of generality, we assume $b>a$. By Corollary \ref{thmcone}, it follows that $\theta_1(t)\in [a,\, b]$ and $\theta_2(t)\in [a+\frac{\pi}{2},\,b++\frac{\pi}{2}]$. Note that $\theta_i(0)$ is a minimum of $\theta_i(t) \, (i=1,2)$ for $t \in [0,\, 1]$.  It implies that $\dot{\theta}_1(0) \geq 0$ and $\dot{\theta}_2(0) \geq 0$. If $\dot{\theta}_1(0)=0$ or $\dot{\theta}_2(0)=0$, by $(d)$ of Lemma \ref{lemmatheta}, there exists some small $\delta>0$, such that $\dot{\theta}_1(t) \dot{\theta}_2(t) < 0$ for all $t \in (0, \delta)$. Contradiction! Hence, both $\dot{\theta}_1(0)$ and $\dot{\theta}_2(0)$ are positive. Note that the angular momentum $A_{m}$ satisfies
\begin{eqnarray*}
&A_{m}(t)&= \sum_{i=1}^3 q_i(t) \times \dot{q}_i(t)\\
&=& \frac{1}{4} Z_1(t) \times \dot{Z}_1(t)+ \frac{4m^2+1}{(2m+1)^2} Z_2(t) \times \dot{Z}_2(t)\\
&=&\left(0,\quad 0, \quad  \frac{1}{4} |Z_1(t)|^2 \dot{\theta}_1(t)+  \frac{4m^2+1}{(2m+1)^2}  |Z_2(t)|^2 \dot{\theta}_2(t) \right).
\end{eqnarray*}
It follows that the angular momentum $A_{m}=A_{m}(0)$ is nonzero.
	
Next, we show $\theta_i \, (i=1,2)$ has no local extrema. If not, we assume $\theta_1$ has a local extrema. Note that $\theta_1(0)=a$, $\theta_1(1)=b$, and $\theta_1(t) \in [a, \, b]$ for all $t \in [0,\,1]$. It is clear that $\theta_1$ must have at least two local extrema. Contradiction to Theorem \ref{geometryangle}!  Similarly, we can show $\theta_2$ has no local extrema. Therefore, both $\theta_1$ and $\theta_2$ are monotone.
\end{proof}

\section{Applications to two sets of periodic orbits}\label{applications}
In this section, we study two sets of periodic orbits, which can be characterized as minimizers connecting collinear configurations and isosceles configurations. By applying Theorem \ref{maingeometricresult1} and Theorem \ref{maingeometricresult2}, we can show their variational existence and prove some geometric properties.

Let $Q_S$ be the set of configurations on the $x$-axis with order constraints $q_{2x}(0) \leq q_{1x}(0) \leq q_{3x}(0)$, and $Q_E$ be the set of isosceles triangles, where the symmetry axis of an isosceles is a counterclockwise $\theta$ rotation of the $x$-axis and $q_1$ is its vertex. (See Fig. \ref{picQsQe}.)

For each $m>0$ and $\theta \in [0, \, \pi/2)$, standard variational results \cite{CH3, FT, Yan1, Yan2, Yan3} imply that there exists a minimizer $\mathcal{P}_{m, \, \theta} \in H^1([0,1], \Sigma)$, such that
\begin{equation}\label{actionminsetting0}
\mathcal{A}(\mathcal{P}_{m, \, \theta})= \inf_{q \in  P(Q_S,\, Q_E)} \mathcal{A}(q) = \inf_{q \in  P(Q_S, \, Q_E)} \int_0^1 (K+U) \, dt,
\end{equation}
where
\[ \mathrm P(Q_S,\, Q_E)=\left  \{q\in H^1([0,1],\, \Sigma) \, \big{|}  \, q(0) \in Q_S, \, q(1) \in Q_E \right \}.\]
The minimizer $\mathcal{P}_{m, \, \theta}$ could have collision singularities. In fact, it is known \cite{Mar, CA, CH3, Ven, TV, Ven2, Yan3} that $\mathcal{P}_{m, \, \theta}$ has no total collision and it is collision-free in $(0,1]$. However, due to the order constraint $q_{2x}(0) \le q_{1x}(0) \le q_{3x}(0)$ in $Q_S$, it is not easy to exclude the possible binary collisions at $t=0$. In the following proposition, we apply Theorem \ref{maingeometricresult1} and Theorem \ref{maingeometricresult2} to show that $\mathcal{P}_{m, \, \theta}$ is collision-free and has some interesting properties.

\begin{proposition}\label{mainthm1}
For each given $m>0$ and $\theta \in [0, \pi/2)$, a minimizer $\mathcal{P}_{m, \, \theta}$ in \eqref{actionminsetting0} is collision-free and it can be extended to a periodic or quasi-periodic orbit. Furthermore, if $\mathcal{P}_{m, \, \theta}$ does not coincide with an Euler orbit, then the polar angles $\theta_1(t)$ and $\theta_2(t)$ of the Jacobi coordinates in $\mathcal{P}_{m, \, \theta}$ have at most one critical point.
\end{proposition}
\begin{proof}
By previous works \cite{Mar, CA, CH3, Ven, TV, Ven2, Yan3} on collision singularities in minimizers, it is known that $\mathcal{P}_{m, \, \theta}$ has no total collision and it is collision-free in $(0, \, 1]$. We are left to exclude possible binary collisions at $t=0$.

Note that at $t=0$, the three masses are on the $x$-axis with an order constraint $q_{2x}(0) \le q_{1x}(0) \le q_{3x}(0)$. It implies that the only possible binary collisions are $q_1(0)=q_2(0)$ and $q_1(0)=q_3(0)$. Since $\mathcal{P}_{m, \, \theta}$ is free of total collision, it follows that $q_{2x}(0)<0$ and $q_{3x}(0)>0$.

If $Z_{2x}(0) = (1+\frac{1}{2m})q_{1x}(0)=0$, it is an Euler configuration. In this case, a binary collision at $t=0$ implies a total collision. Contradiction! It follows that $\mathcal{P}_{m, \, \theta}$ has no binary collision at $t=0$ when $Z_{2x}(0)=0$.

 If $Z_{2x}(0) > 0$, the only possible binary collision is $q_1(0)=q_3(0)$. Similarly, if $Z_{2x}(0) = (1+\frac{1}{2m})q_{1x}(0) < 0$,  the only possible binary collision is $q_1(0)=q_2(0)$. Here we only discuss the case when $Z_{2x}(0)> 0$. The other case $Z_{2x}(0)< 0$ follows by a similar argument.

We prove it by contradiction. Assume $q_1(0)=q_3(0)$ in $\mathcal{P}_{m, \, \theta}$. We can analyze the asymptotic behavior of the minimizer at $t=0$. In fact, by \cite{CH3, FU, Yu1, Yan3}, the following limit holds:
\begin{equation}\label{yvelocity}
 \lim_{t \to 0^{+}} \dot{q}_{1y}=\lim_{t \to 0^{+}} \dot{q}_{3y}= -\frac{m+1}{m} \dot{q}_{2y}(0).
\end{equation}
When $\dot{q}_{2y}(0)= 0$, by Lemma 6.2 in \cite{Yu1}, the collision path $Z=(Z_1, Z_2) = (Z_1(t), Z_2(t))(t\in [0,1])$ must stay on the $x$-axis with order $q_{2x}(t) \leq q_{3x}(t) \leq q_{1x}(t)$. However, by the definition of $Q_E$, the collision-free isosceles configuration at t = 1 can never become a collinear configuration on the x-axis with such an order. Contradiction!\\
 When $\dot{q}_{2y}(0)>0$, By \eqref{yvelocity},  it follows that $\displaystyle \lim_{t \to 0^{+}} \dot{q}_{1y} < 0$. Then for $t\in (0,\epsilon]$ with $\epsilon>0$ sufficiently small, we have
 \[\dot{q}_{1y}(t)<0, \qquad \forall t \in (0, \epsilon]. \]
 Hence there exists some $\epsilon>0$ small enough, such that $Z_1(t) \in \mathsf{Q}_4$ and $Z_2(t) \in \mathsf{Q}_4$ for $t \in (0, \, \epsilon)$. A similar conclusion holds when $\dot{q}_{2y}(0)<0$.

 On the other hand, Theorem \ref{maingeometricresult1} implies that $Z_1(t)$ and $Z_2(t)$ are always in two adjacent quadrants for $t \in [0, \, 1]$. We prove it in two cases.

When $\theta =0$, we have $Z_1(1)$ is on the $y$-axis and $Z_2(1)$ is on the $x$-axis.  By Theorem \ref{maingeometricresult1}, it follows that
\begin{equation*}
Z_1(t)\in \overline{\mathsf{Q}_4},\ \ \ Z_2(t)\in \overline{\mathsf{Q}_1},\ \ \ \forall t\in [0,1],
\end{equation*}
or
\begin{equation*}
Z_1(t)\in \overline{\mathsf{Q}_1},\ \ \ Z_2(t)\in \overline{\mathsf{Q}_4},\ \ \ \forall t\in [0,1].
\end{equation*}

When $\theta\in (0,\pi/2)$, $Z_1(1)\in \mathsf{Q}_2 \, \,  \text{or} \, \, \mathsf{Q}_4$ and $Z_2(1) \in \overline{\mathsf{Q}_1} \, \, \text{or} \, \,  \overline{\mathsf{Q}_3}$. By Theorem \ref{maingeometricresult1}, it implies that
\[Z_1(t) \in  \overline{\mathsf{Q}_4}, \qquad   Z_2(t) \in \overline{\mathsf{Q}_1}. \]
Hence, for any $m>0$ and $\theta \in [0, \, \pi/2)$, $Z_1(t)$ and $Z_2(t)$ are always in two adjacent quadrants for $t \in [0, \, 1]$. However, the asymptotic behavior at $t=0$ implies that $Z_1(t)$ and $Z_2(t)$ are in the same quadrant for $t \in (0, \, \epsilon)$. Contradiction!

Therefore, $\mathcal{P}_{m, \, \theta}$ is free of collision and it is a solution of the Newtonian equations \eqref{newton}. By applying the first variation formula as in Section 5 of \cite{Yan1}, one can show that $\mathcal{P}_{m, \, \theta}$ can be extended to a periodic or quasi-periodic orbit.

In the end, it is clear that $Z=(Z_1,\, Z_2)$ can not be an isosceles motion or a collinear motion. However, it could be part of an Euler orbit. By Theorem \ref{maingeometricresult2}, if $\mathcal{P}_{m, \, \theta}$ does not coincide with an Euler orbit, the polar angles $\theta_1(t)$ and $\theta_2(t)$ of the Jacobi coordinates in $\mathcal{P}_{m, \, \theta}$ have at most one critical point.
\end{proof}

\begin{figure}[ht]
    \begin{center}
    \subfigure[ \, Periodic orbit extended by $\mathcal{P}_{1, \, 0}$   ]{\includegraphics[width=2.5in]{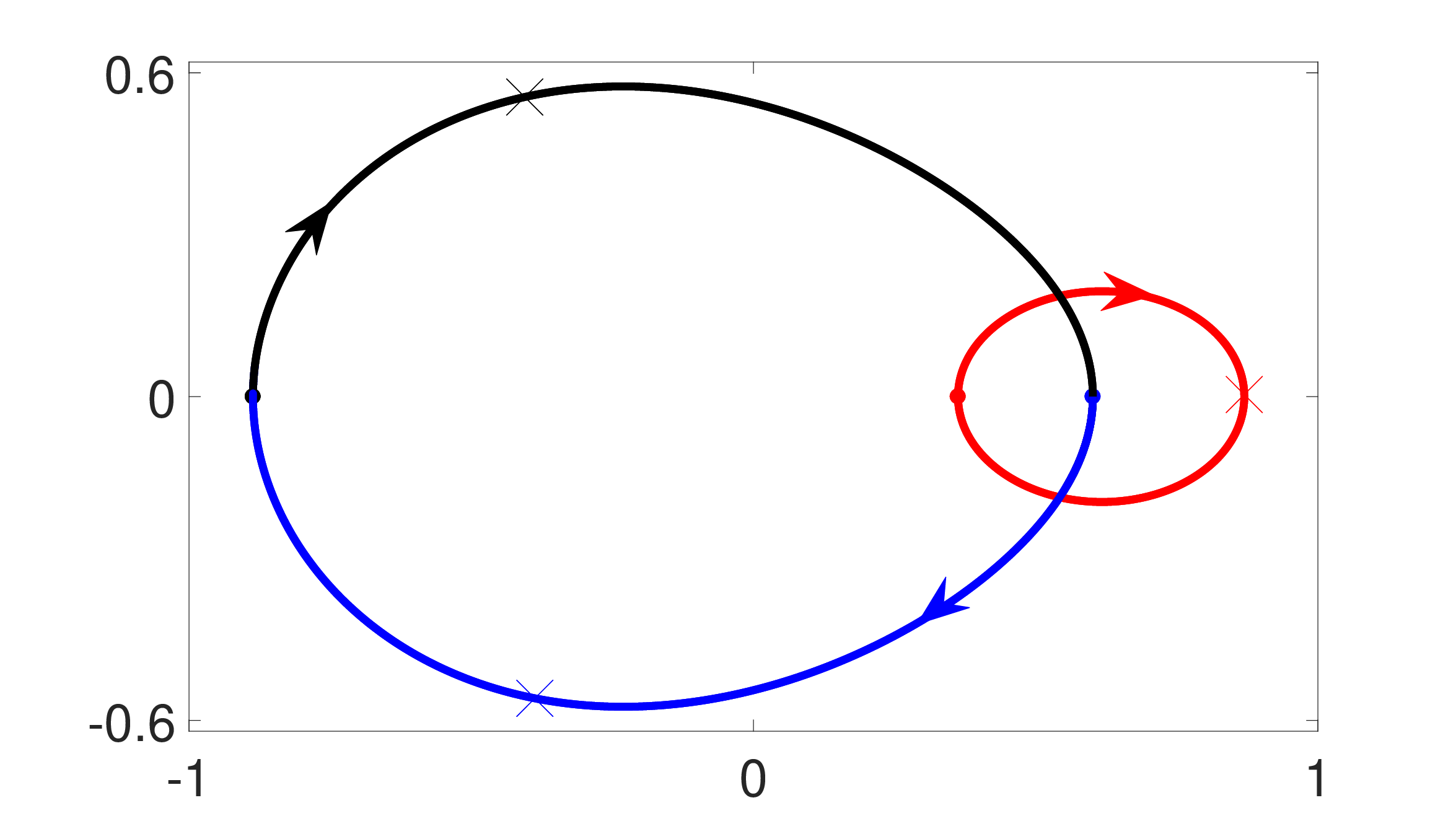}}
    \subfigure[ \, $\theta_2(t)$ has a critical point in $\mathcal{P}_{1, \, 0}$ ]{\includegraphics[width=2.5in]{Z2turn.eps}}
   \end{center}
 \caption{ \label{min1}The picture on the left is the periodic orbit extended by the minimizer $\mathcal{P}_{1, \, 0}$, while its right is $\mathcal{P}_{1, \, 0}$ under the Jacobi coordinates.}
 \end{figure}

The following figure shows the graphs of a minimizer $\mathcal{P}_{1, \, 0}$ under both the Cartesian coordinates and the Jacobi coordinates. Fig. \ref{min1} (a) is $\mathcal{P}_{1, \, 0}$ under the Cartesian coordinates and its periodic extension. Fig. \ref{min1} (b) is $\mathcal{P}_{1, \, 0}$ under the Jacobi coordinates. Numerical result implies that $\theta_1(t)$ has no critical point, while $\theta_2(t)$  has exactly one critical point in $\mathcal{P}_{1, \, 0}$.

If we change the order in the collinear configurations of $Q_S$, it leads to a different set of periodic orbits. Let $Q_{S_1}$ be the set of collinear configurations on the $x$-axis with order constraints $q_{1x}(0) \leq q_{2x}(0) \leq q_{3x}(0)$. For each $\theta \in (0, \pi/2]$, there exists a minimizer $\widetilde{\mathcal{P}_{m, \, \theta}} \in H^1([0,1], \Sigma)$, such that
\begin{equation}\label{actionminsetting12}
\mathcal{A}(\widetilde{\mathcal{P}_{m, \, \theta}})= \inf_{q \in  P(Q_{S_1}, \, Q_{E})} \mathcal{A}(q) = \inf_{q \in  P(Q_{S_1}, \, Q_{E})} \int_0^1 (K+U) \, dt.
\end{equation}

Similarly, we can show that
\begin{proposition}\label{retrograderesult}
For each given $m>0$ and $\theta \in (0, \, \pi/2)$, a minimizer $\widetilde{\mathcal{P}_{m, \, \theta}}$ in \eqref{actionminsetting12} is collision-free, and it can be extended to a periodic or quasi-periodic orbit. Furthermore, the Jacobi coordinates $Z=(Z_1, \, Z_2)$ of $\widetilde{\mathcal{P}_{m, \, \theta}}$ satisfy that $Z_1 \in \overline{\mathsf{Q}_4}$ and $Z_2 \in \overline{\mathsf{Q}_3}$, which indicates that $\widetilde{\mathcal{P}_{m, \, \theta}}$ contains no collinear configuration except for the boundaries. The corresponding polar angles $\theta_1(t)$ and $\theta_2(t)$ have at most one critical point.
\end{proposition}

For $\theta \in (0, \pi/2)$, the orbit generated by $\widetilde{\mathcal{P}_{m, \, \theta}}$ is similar to the retrograde orbit in \cite{CH1}. Mathematically, they may not be the same since the orbit in Proposition \ref{retrograderesult} has stronger symmetry than the retrograde orbit in  \cite{CH1}. Indeed, by the first variation formula, $\widetilde{\mathcal{P}_{m, \, \theta}}$ has velocities perpendicular to $x$-axis at $t=0$. However, in the settings of \cite{CH1}, it does not have to be so.

When $\theta=\pi/2$, it is closely related to one of the open problems \cite{VEOP} proposed by Venturelli in 2003. We have shown in \cite{Yan3} that $\widetilde{\mathcal{P}_{m, \, \pi/2}}$ coincide with either the Schubart orbit or the Broucke-H\'{e}non orbit. It will be interesting if one can show that  $\widetilde{\mathcal{P}_{m, \, \pi/2}}$ is collision-free. Numerically, Venturelli \cite{VEOP} claimed that $\widetilde{\mathcal{P}_{m, \, \pi/2}}$ is collision-free at least for $m=1$.

\begin{figure}[ht]
    \begin{center}
    \subfigure[ \, Broucke-H\'{e}non orbit extended by $\widetilde{\mathcal{P}_{1, \, \pi/2}}$   ]{\includegraphics[width=2.5in]{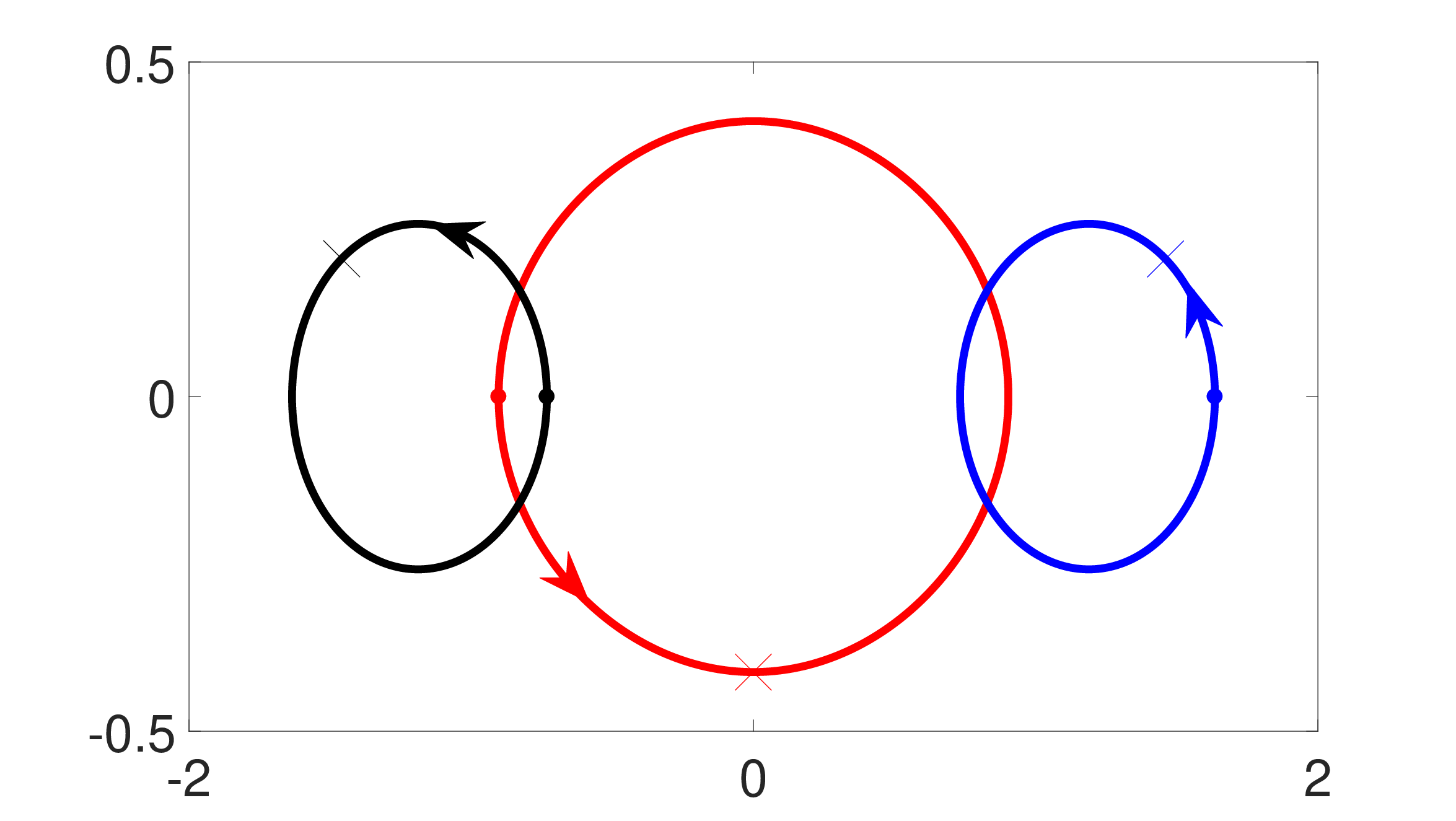}}
    \subfigure[ \, $\theta_1(t)$ has a critical point in $\widetilde{\mathcal{P}_{1, \, \pi/2}}$]{\includegraphics[width=2.5in]{Z1turn.eps}}
        \subfigure[ \, Periodic orbit extended by $\widetilde{\mathcal{P}_{1, \, \pi/6}}$   ]{\includegraphics[width=2.5in]{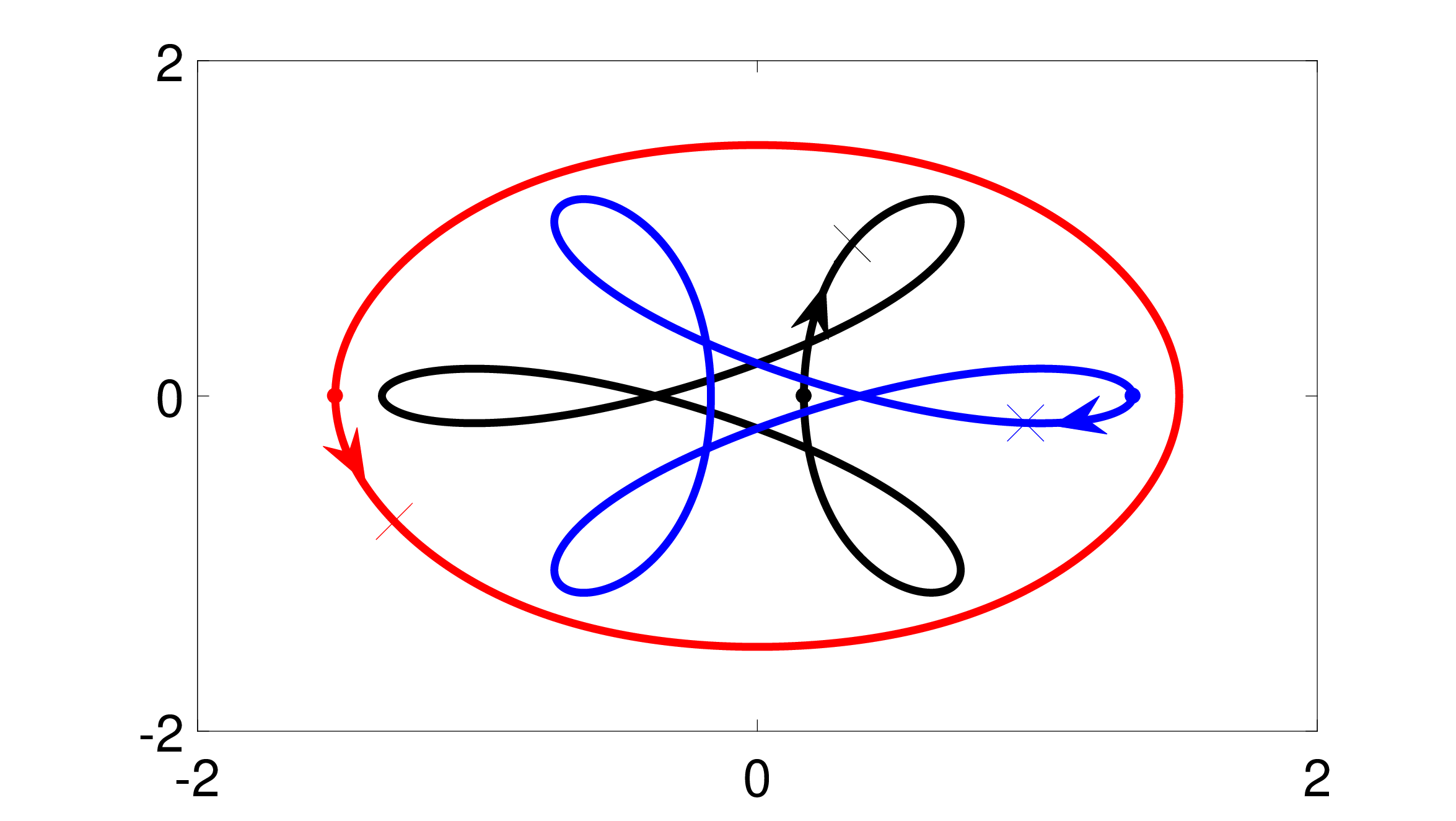}}
    \subfigure[ \, $\theta_1(t)$, $\theta_2(t)$ have no critical point in $\widetilde{\mathcal{P}_{1, \, \pi/6}}$]{\includegraphics[width=2.5in,height=1.2in]{Z1Z2noturn.eps}}
   \end{center}
 \caption{ \label{min2}The left two graphs are periodic orbits extended by the two minimizers $\widetilde{\mathcal{P}_{1, \, \pi/2}}$ and $\widetilde{\mathcal{P}_{1, \, \pi/6}}$, while the right ones are the two minimizers under the Jacobi coordinates.}
 \end{figure}

Similar to Fig. \ref{min1}, we can draw the pictures of two minimizers $\widetilde{\mathcal{P}_{1, \, \pi/2}}$ and $\widetilde{\mathcal{P}_{1, \, \pi/6}}$ under both the Cartesian coordinates and the Jacobi coordinates in Fig. \ref{min2}. On its left, (a) and (c) are the graphs of the two minimizers under the Cartesian coordinates and their periodic extensions. On its right, (b) and (d) are the two minimizers under the Jacobi coordinates. In (b), $\theta_1(t)$ has one critical point. While in (d), both $\theta_1(t)$ and $\theta_2(t)$ have no critical point. The numerical results in Fig. \ref{min1} and Fig. \ref{min2} indicate that $\theta_i(t) \, (i=1,2)$ could have no critical point or $\theta_i(t) \, (i=1, \, \text{or} \, 2)$ has one critical point. Therefore, numerical investigation suggests that Theorem \ref{maingeometricresult2} is a sharp result.

\section*{Acknowledgements}
The authors gratefully acknowledge the support of NSFC (No. 11901279, 11871086). We thank all the referees for their time and efforts! 

\end{document}